\documentclass[review]{elsarticle}

\usepackage{amsmath}
\usepackage{latexsym, amssymb}
\usepackage{graphicx}
\usepackage{amsmath, amsbsy}
\usepackage{amsopn, amstext}
\usepackage{hyperref}
\usepackage{cancel, color}
\usepackage{epstopdf}

\def\J{{\bf 1}}

\DeclareMathOperator{\rank}{rank}

\DeclareMathOperator{\Span}{Span}
\DeclareMathOperator{\Col}{Col}

\DeclareMathOperator{\lcm}{lcm}

\def\cal{\mathcal}

\def\ra{\rightarrow}

\def\a{\alpha}
\def\b{\beta}
\def\d{\delta}

\def\D{\Delta}

\def\0{{\bf 0}}

\newcommand{\R}{{\mathbb R}}
\newcommand{\C}{{\mathbb C}}
\newcommand{\N}{{\mathbb N}}

\newcommand{\F}{{\mathbb F}}

\def\dsum{\mathop{\sum}\limits}

\newtheorem{thm}{Theorem}[section]
\newtheorem{dfn}[thm]{Definition}
\newtheorem{prp}[thm]{Proposition}
\newtheorem{exa}[thm]{Example}

\newtheorem{cor}[thm]{Corollary}
\newtheorem{rem}[thm]{Remark}
\newtheorem{alg}[thm]{Algorithm}

\begin{document}

\begin{frontmatter}

\title{On Universal Eigenvalues and Eigenvectors  of Hypermatrices\tnoteref{footnoteinfo}}

\tnotetext[footnoteinfo]{ This work is supported partly by the National Natural Science Foundation of China (NSFC) under Grant 62073315 and 62350037. Corresponding author: Zhengping Ji.}

\author[AMSS,STPC]{Daizhan Cheng}\ead{dcheng@iss.ac.cn}
\author[AMSS,CAS]{Zhengping Ji}\ead{jizhengping@amss.ac.cn}

\address[AMSS]{Key Laboratory of Systems and Control, Academy of Mathematics and Systems Science, Chinese Academy of Sciences, Beijing 100190, P.R.China}
\address[STPC]{Research Center of Semi-tensor Product of Matrices, Theory and Applications, Liaocheng University, Liaocheng, P.R. China}
\address[CAS]{School of Mathematical Sciences, University of Chinese Academy of Sciences, Beijing 100049, P.R.China}

\begin{abstract}
A generalized eigenvector of a hypermatrix, called the universal (U-) eigenvector, is proposed, which extended the notion of diagonal (D-) eigenvectors in the literature.
Using the semi-tensor product, the homogeneous U-eigenequation can be converted into a general eigenequation of matrix  $(A-\lambda B)x=0$. A general technique for solving this equation is proposed, which leads to two kinds of eigenvalues: essential and quasi eigenvalues. The technique to convert nonhomogeneous eigenequation to homogeneous ones is also revealed.
Then a hypervector decomposing method, called the monic decomposition algorithm (MDA), is developed.
Using the MDA,  the U-eigenproblem (including the D-eigenproblem) can be converted into general matrix eigenproblems.
Some examples are presented, demonstrating the geometric meaning and potential applications of the U-eigenvalue/eigenvector.
\end{abstract}

\begin{keyword}
	Hypermatrix, Hypervector, U- (D-) eigenproblem, essential (quasi-) eigenvalue,  Semi-tensor product of matrices.
\end{keyword}

\end{frontmatter}

\section{Introduction}

The eigenvalues and eigenvectors of matrices play important roles in theoretical analysis and applications of matrices. A matrix (or a vector) is essentially a set of order $2$ data (respectively, a vector is of order $1$). When a set of ordered data with order $d>2$ is considered, it becomes a hypermatrix \cite{lim13}. Extending the notions of eigenvalues and eigenvectors to the hypermatrix setting becomes natural. In 2005, Qi \cite{qi05} and Lim \cite{lim05} firstly introduced the definition of tensor eigenvalue/eigenvector independently. Since then, studies has been devoted to the eigenvalue/eigenvector problems of hypermatrices (or tensors).

For different purposes, various kinds of eigenvalues/eigenvectors of hypermatrices have been developed. We list some commonly used definitions of the eigenvalues/eigenvectors  as follows. Let $A=\{a_{i_1,\cdots,i_d}\in \R\;|\;i_s\in [1,n], s\in [1,d]\}$ be a  hypermatrix \cite{lim13}, and $x=(x_1,\cdots,x_n)\in \R^n$ (or $\C^n$). In \cite{qi05} the author defined a map as
\begin{align}\label{1.1.1}
	[M_Ax^{d-1}]_{i}:=\dsum_{i_2=1}^n\cdots\dsum_{i_d=1}^na_{i,i_2,\cdots,i_d}x_{i_2}\cdots x_{i_d},\quad i\in [1,n].
\end{align}
Here $M_Ax^{d-1}\in \R^n$ is considered as a column vector. Several types of eigenvalues/eigenvectors have been defined for $M_A$ as follows.

\begin{itemize}
	\item[(i)] The H-eigenvalue and eigenvector \cite{qi05}:
	\begin{align}\label{1.1.2}
		M_Ax^{d-1}=\lambda x^{[d-1]}, \quad  x^{[d-1]}:=(x_1^{d-1},\cdots,x_n^{d-1})^T.
	\end{align}
	\item[(ii)] The E- (or Z-) eigenvalue and eigenvector \cite{qi05}:
	\begin{align}\label{1.1.3}
		M_Ax^{d-1}=\lambda x ~\mbox{with}~x^Tx=1.
	\end{align}
	\item[(iii)] The D-eigenvalue and eigenvector \cite{qi08}:
	\begin{align}\label{1.1.4}
		M_Ax^{d-1}=\lambda Dx ~\mbox{with}~x^TDx=1.
	\end{align}
	\item[(iv)] The Markov chain-based eigenvalue and eigenvector \cite{chi13,li14}:
	\begin{align}\label{1.1.5}
		M_Ax^{d-1}=\lambda (x_1+\cdots+x_n)^{d-2}x.
	\end{align}
	\item[(v)] The generalized  eigenvalue and eigenvector \cite{din15}:
	\begin{align}\label{1.1.6}
		M_Ax^{d-1}=\lambda Bx^{d-1}.
	\end{align}
		
\item[(vi)] The cubic eigenvalue and eigenvector with inner product \cite{kol14}:
\begin{align}\label{1.1.7}
M_Ax^3=\lambda (x^Tx)x.
\end{align}

\end{itemize}

The different types of eigenproblems  of hypermatrices have been explored and widely applied to practical problems. For instance, the  properties of tensors via eigenvalue/eigenvector approach has been studied in \cite{cha08}, \cite{cha09}. In particular, some   new properties for non-negative tensors were revealed in \cite{cha13}. \cite{cui14} considered the symmetric tensors and provided two numerical algorithms for generalized tensor eigenproblems. \cite{kol14} proposed an adaptive shifted power method to convert the solvability of eigenproblems to solving an optimization problems.  \cite{ni14} proposed the US-eigenvalue, which is raised from quantum information processing. As a special case of more general definition, it seems more suitable for quantum mechanical systems. Higher-order singular value decomposition \cite{del00} and tensor singular value decomposition \cite{pan22} were developed for the eigenvalue/eigenvector based hypermatrix analysis. etc., just to mention a few.

Recently, the development of dynamic control systems over order $3$ hypermatrices is emerging, and has given rise to many practical uses in physics \cite{chan22}, image processing \cite{arz18,kil13}, and biological systems \cite{yu22}, etc. Similar to the role played by the eigenproblems  in systems governed by linear ordinary differential equations, the eigenvalues/eigenvectors of hypermatrices is a crucial concept in the analysis of dynamic control systems over hypermatrices \cite{cche24}.

The semi-tensor product (STP) is the main tool in our new approach.  The STP of matrices is a generalization of the classical matrix product to two arbitrary matrices.  STP has been proposed two decades before \cite{che01}.  Since then, it has been rapidly developed in both theoretical aspect and various applications \cite{che11,che12}. For instance, it has been applied to studying Boolean networks and finite valued networks (see survey papers \cite{for16,li18,lu17,muh16}); finite games (see the survey paper \cite{che21}); finite automata (see survey paper \cite{yan22}); dimension-varying systems \cite{che19,che19b}, etc.

{\bf Throughout this paper the default matrix product is assumed to be semi-tensor product (STP)}, denoted by $\ltimes$, unless elsewhere is stated. STP  will be discussed in the next section. For the reader, who is not familiar with it,  all needed so far is the following fact, which is also important in the sequel.

\begin{prp}\label{p1.1.0}
Let  $x$ and $y$ be two column vectors. Then
\begin{align}\label{1.1.701}
x\ltimes y=x\otimes y.
\end{align}
\end{prp}
This paper propose a new kind of eigenvector, called U-eigenvector, which is a hypervector, $x=\ltimes_{i=1}^rx_i$, where $x_i\in \R^n$, $i\in[1,r]$. (In column vector case, $x\ltimes y=x\otimes y$.) A process is proposed to solve U-eigenproblem, which is described by Figure \ref{Fig.0.1}.

\begin{figure}
\centering
\setlength{\unitlength}{0.8 cm}
\begin{picture}(16,14)(0,-2)
\thinlines
\put(0,-2){\framebox(8,1){$x_i, ~i\in [1,r]$: Solution of (\ref{0.0.1})}}
\put(0,3){\framebox(8,1){$x$ : Solution of (\ref{0.0.2})}}
\put(0,6){\framebox(8,1){$z$: Solution of (\ref{0.0.4})}}
\put(0,9){\framebox(8,1){General Matrix Equiequation (\ref{0.0.1})}}
\put(0,12){\framebox(8,1){D-Equiequation (\ref{0.0.2})}}
\put(0,15){\framebox(8,1){U-Equiequation (\ref{0.0.4})}}
\thicklines
\put(4,15.2){\vector(0,-1){2.2}}
\put(4,12.2){\vector(0,-1){2.2}}
\put(4,9.2){\vector(0,-1){2.2}}
\put(4,6.2){\vector(0,-1){2.2}}
\put(4,3.2){\vector(0,-1){4.2}}
\put(4.2,13.8) {STP Transfer}
\put(4.2,10.8) {Raising Power Algorithm}
\put(4.2,7.8) {Solving general Matrix Eigenproblem}
\put(4.2,4.8) {Lowering Power Algorithm}
\put(4.2,1.8) {Monic Decomposition Algorithm}
\put(4.2,0.8) {$+$}
\put(4.2,-0.2) {Iteration Algorithm}
\end{picture}
\caption{Solving Process\label{Fig.0.1}}
\end{figure}

\begin{itemize}
\item[]  U-Eigenequation

\begin{align}\label{0.0.1}
\begin{cases}
A\ltimes_{i=1}^rx_i=B\ltimes_{j=1}^sy_j,\\
y_j=B_j\ltimes_{i=1}^rx_i.
\end{cases}
\end{align}

\item[]  D-Eigenequation

\begin{align}\label{0.0.2}
Ax^r=\lambda Bx^s.
\end{align}

Particularly, the D-eigenequation obtained from (\ref{0.0.1}) is as follows.
\begin{align}\label{0.0.3}
Ax=\lambda Bx^s,\quad x=\ltimes_{i=1}^rx_i.
\end{align}

\item[] General Matrix Eigenequation\cite{van75}

\begin{align}\label{0.0.4}
Az=\lambda Bz.
\end{align}

Particularly, the general matrix eigenequation obtained from (\ref{0.0.2}) is with $z=x^{t}$, where $t=\max(r,s)$.
\end{itemize}

The main work of this paper consists of the following:
\begin{itemize}
\item[(i)] Using STP,  the  STP transfer, the Raising and Dowering Power Algorithm are proposed.
\item[(ii)] A Monic Decomposition Algorithm has been developed.
\end{itemize}
The above two kinds of algorithms are straightforward computable. Hence, the {\bf U-eigenproblem is converted to general matrix eigenproblem}. Two kinds of eigenvalues, called essential and quasi eigenvalues, are proposed, and the algorithm has also been investigated. Since the general matrix eigenproblem is of a matrx equetian

It is obvious that all the eigenequations in literature so far are all D-eigenproblems. So the general solving process can also be used to solve  them.

The rest of this paper is organized as follows. Section 2 formulates the U-eigenproblem. Section 3 provides some mathematical preliminaries, including  (i) a brief review of the STP (ii) matrix expression of hypermatrix; (iii) the hypervector, which is the form of U-eigenvector, is introduced. In Section 4 the
type ${\cal B}$ and general U-eigenproblem is discussed first. Then the D-eigenproblem has also been investigated.  For each case the linear algebraic expression is constructed, which means to turn the eigenequations into the general matrix eigenequation.   Section 5 provides the monic decomposition algorithm. Using it,  we can verify whether the $x$ obtained in Section 4 is a hypervector. In this way, we can find the U-or  D- eigenvalue-eigenvectors. In Section 6 some illustrative examples are presented. first example demonstrates the geometric meaning of U-eigenvectors. Second example comes from \cite{kol14,din15}. It shows the advantage of our approach over exiting method.
Third one shows its application to tensor. Last two are  neural network related. Last three examples  show the potential applications of U-eigenvalues and eigenvectors.
  Section 7 is a brief conclusion with two remarks. (i) The difference of hypermatrix with tensor.  (ii) Relationship between eigenproblem of hypermatrix and the non-square matrix theory.

Before ending this section, we give a list of notations used in this paper.

\begin{enumerate}
\item $\N$: Set of positive integers.
\item $\F^n$: $n$ dimensional Euclidean space over $\F$, where $\F=\R$ or $\F=\C$.
\item ${\cal M}_{m\times n}$: the set of $m\times n$ matrices (over $\F$).
\item $\lcm(a,b)$: least common multiple of $a$ and $b$, ($a,b\in \N$).
\item  $\F^{n_1\times n_2\times \cdots \times n_d}$: the $d$-th order hypermatrices of dimensions $n_1,n_2,\cdots,n_d$.
\item  $\F^{n_1\ltimes n_2\ltimes \cdots \ltimes n_d}$: the $d$-th order hypervectors of dimensions $n_1,n_2,\cdots,n_d$.
\item $[a,b]$: the set of integers $a\leq i \leq b$.
\item $\d_n^i$: the $i$-th column of the identity matrix $I_n$.
\item $\d_n^I:=[(\d_n^1)^\mathrm{T}, (\d_n^2)^{\mathrm{T}},\cdots,(\d_n^n)^\mathrm{T}]^\mathrm{T}$.
\item $\D_n:=\left\{\d_n^i\vert i=1,\cdots,n\right\}$.
\item $\d_n[i_1,\cdots,i_s]:=\left[\d_n^{i_1},\cdots,\d_n^{i_s}\right]$.
\item $\J_{\ell}:=(\underbrace{1,1,\cdots,1}_{\ell})^\mathrm{T}$.
%
%
\item $S(x)$: the subspace spanned by $\{x_1,\cdots,x_r\}$, where $x=\ltimes_{i=1}^rx_i$.
\item $\otimes$:  Kronecker product of matrices.
\item $\ltimes$: semi-tensor product of matrices.
\item $\times_{{\bf j}}$: contraction product of two hypermatrices with respect to  their common index subset ${{\bf j}}$.
\item ${\cal T}^r_s(V)$: the set of tensors over $V$ with covariant order $r$ and contra-variant order $s$. When $s=0$,
${\cal T}^r(V)$ denotes the set of covariant tensor over $V$ with covariant order $r$.
\end{enumerate}

\section{Problem Formulation}

The  U-eigenvector, with respect to (w.r.t.) an eigenvalue proposed in this paper  is a hypervector, which is expressed as
$$
x=\ltimes_{i=1}^rx_i,\quad x_i\in \R^n.
$$
When $x_i=z$, $i\in [1,r]$, $x$ is called a D- eigenvector. Note that if $x=z^r$ is an U-eigenvector, then $z$ is the ``eigenvector" used so far in literature. But in this paper $z^r$ is called the eigenvector.

The first motivation for U-eigenproblem comes from  the observation on the contraction product of hypermatrices.

\begin{dfn}\label{d1.1.1} \cite{lim13}
\begin{itemize}
\item[(i)] A hypermatrix of order $d$ and dimension $n_1\times \cdots n_d$ is defined as
$$
{\cal A}=\{a_{i_1,i_2,\cdots,i_d}\;|\; i_j\in[1,n_j],\; j\in [1,d]\}\in \R^{n_1\times \cdots \times n_d}.
$$
When $n_i=n$, $i\in [1,d]$, ${\cal A}$ is called an equilateral hypermatrix.

\item[(ii)] Assume ${\cal A}=(a_{i_1,i_2,\cdots,i_r})\in \R^{m_1\times \cdots \times m_r}$,
${\cal B}=(b_{j_1,j_2,\cdots,j_s})\in \R^{n_1\times \cdots \times n_s}$. Denote the index sets as
$$
{\bf i}=\{i_1,i_2,\cdots,i_r\};\quad {\bf j}=\{j_1,j_2,\cdots,j_s\},
$$
and
$$
{\bf k}=\{k_1,k_2,\cdots, k_t\}\subset {\bf i}\bigcap {\bf j}.
$$
Moreover, assume
$$
m_{i_{k_{\ell}}}=n_{j_{k_{\ell}}},\quad \ell\in [1,t].
$$
For the sake of notational ease, we rename the index sets as
$$
\begin{array}{l}
{\bf i}=\{i_1,\cdots,i_p,k_1,\cdots,k_t\},\quad p=r-t,\\
{\bf j}=\{j_1,\cdots,j_q,k_1,\cdots,k_t\},\quad q=s-t.\\
\end{array}
$$
Then the contraction product of ${\cal A}$ and ${\cal B}$ with respect to ${\bf k}$ is defined as follows.
\begin{align}\label{1.1.8}
{\cal A}\times_{{\bf k}}{\cal B}={\cal C}\in \F^{m_{i_1}\times \cdots\times m_{i_p}\times n_{j_1}\times \cdots\times n_{j_q}},
\end{align}
where
$$
c_{i_1,\cdots,i_p,j_1,\cdots,j_q}=\dsum_{k_1=1}^{m_{k_1}}\cdots \dsum_{k_t=1}^{m_{k_t}}
a_{i_1,\cdots,i_p,k_1\cdots k_t} b_{j_1,\cdots,j_q,k_1\cdots k_t}.
$$
\end{itemize}
\end{dfn}	

Assume ${\cal A}=(a_{i_1,\cdots,i_d}) \in \F^{\overbrace{n\times \cdots\times n}^d}$,
$x=(x_{j_1,\cdots,j_r}) \in \F^{\overbrace{n\times \cdots\times n}^r}$, $1\leq r\leq d-1$. Denote
$$
{\bf i}=\{i_1,\cdots,i_d\},\quad {\bf j}=\{j_1,\cdots,j_r\}.
$$
Assume  ${\bf j}\subset {\bf i}$ and
$$
{\bf k}=\{k_1,\cdots,k_s\}={\bf i}\backslash {\bf j},
$$
where $s=d-r$.
Then ${\cal A}$ can be considered a mapping from $\F^{\overbrace{n\times \cdots\times n}^r}$ to $\F^{\overbrace{n\times \cdots\times n}^s}$ by
\begin{align}\label{1.1.9}
x\mapsto {\cal A}\times_{\bf j}x, \quad x\in \F^{\overbrace{n\times\cdots\times n}^r}.
\end{align}

Based on the contraction product of hypermatrices we give the  hyperproblem of hypermatrices as follows.

\begin{dfn}\label{d1.1.2} Let ${\cal A}=(a_{i_1,\cdots,i_d}) \in \F^{\overbrace{n\times \cdots\times n}^d}$,
$x=(x_{j_1,\cdots,j_r})\in \F^{\overbrace{n\times \cdots\times n}^r}$, where $d-r=s>0$ and ${\bf j}\subset {\bf i}$.
Then $x\neq 0$ is called an eigenvector of ${\cal A}$  of type ${\cal B}$  w.r.t. eigenvalue $\lambda$, if
\begin{align}\label{1.1.10}
{\cal A}\times_{{\bf j}} x=\lambda {\cal B}(x),
\end{align}
where ${\cal B}:\F^{\overbrace{n\times \cdots\times n}^r}\ra \F^{\overbrace{n\times \cdots\times n}^s}$ is a pre-assigned mapping, called the type of the eigenproblem.
\end{dfn}

We can arrange entries of $x$ into a column vector. Then It is easy to verify that $x$ is an eigenvector of ${\cal A}$  w.r.t. eigenvalue $\lambda$ and type ${\cal B}$, if and only if, $x$ makes Figure \ref{Fig.1.1} commutative.

\begin{figure}
\centering
\setlength{\unitlength}{0.6 cm}
\begin{picture}(10,6)
\thicklines
\put(3.7,0.8){$\F^{\overbrace{n\times \cdots\times n}^{s}}$}
\put(0.2,4.8){$\F^{\overbrace{n\times \cdots\times n}^{r}}$}
\put(7.2,4.8){$\F^{\overbrace{n\times \cdots\times n}^{s}}$}
\put(3.5,5){\vector(1,0){3.5}}
\put(2,4.2){\vector(1,-1){2}}
\put(6,2.2){\vector(1,1){2}}
\put(5,5.2){${\cal A}$}
\put(2,2.5){${\cal B}$}
\put(7.2,2.5){$\lambda$}
\end{picture}
\caption{Eigenvalue-Eigenvector of a Hypermatrix ${\cal A}$\label{Fig.1.1}}
\end{figure}

\begin{rem}\label{r1.1.3}
\begin{itemize}
\item[(i)] All the products of hypermatrices used literature for eigenproblems can be expressed as the contraction product. For instance, the product in (\ref{1.1.1}) can be expressed as
$$
M_Ax^{d-1}=A\times_{i_2,i_3,\cdots,i_d} (\ltimes_{j=2}^d x^j),
$$
where $x^j=x$, $j\in[2,d]$,  and we use  $x^j_{i_j}$ to label the $i_j$ th element of $x^j$, which is the $j-1$-th $x$ in $\underbrace{x\ltimes \cdots \ltimes x}_{d-1}$.

\item[(ii)]
The type is a key factor characterizing the eigenproblem of a hypermatrix. It is a generalization of the $B$ in Definition (\ref{1.1.6}). Hence, the role played by ${\cal B}$ can be found partly in \cite{din15}.
In fact, we do not require the ${\cal B}$ being a linear mapping. It will be specified later.

\item[(iii)] As in \cite{din15}, ${\cal B}$ might be a set of mappings, which have  certain common properties.

\item[(iv)] ${\bf j}\subset {\bf  i}$ is the subset of indexes. In contraction product the ${\bf j}$ is the subset of index ${\bf i}$, which will be eliminated in the resulting hypermatrix. Different subset ${\bf j}$ will cause different eigenvalue-eigenvectors.

\item[(v)] When eigenvector $x$ is considered as a special subspace, we need $x=\ltimes_{i=1}^rx_i$. Then $x$ is used to represent a subspace $S(x)\subset \F^n$.
\end{itemize}
\end{rem}

\begin{rem}\label{r1.1.301}
\begin{itemize}
\item[(i)] It is natural to ask: why $x$ is assumed to be a hypervector? The answer is: a hypervector represents a subspace of $\R^n$. This is similar to the eigenvector for a matrix, which represents a one dimensional subspace. Hence the eigenproblem of hypermatrix has obvious  geometric meaning. Say, we will consider whether or not an U-eigenvector represents an invariant subspace, etc. ?
\item[(ii)] From this perspective, the eigenvectors of hypermatrices in the literature so for, i.e., D-eigenvectors, represent exactly one dimensional invariant subspace.
\end{itemize}
\end{rem}

\begin{exa}\label{e1.1.4} Consider a square matrix $A=(a_{i,j})\in {\cal M}_{n\times n}$. Then $d=2$, $r=s=1$. Set ${\cal B}$ as an identity mapping, then
\begin{itemize}
\item[(i)] Choosing ${\bf j}=\{j\}\subset {\bf i}=\{i,j\}$, we have
that
$x\neq 0$ is an eigenvector of $A$ w.r.t. $\lambda$, if and only if
$$
Ax=\lambda x.
$$
\item[(ii)] Choosing ${\bf j}=\{i\}\subset {\bf i}=\{i,j\}$, we have
that
$x\neq 0$ is a right eigenvector of $A$ w.r.t. $\lambda$, if and only if
$$
x^{\mathrm{T}}A=\lambda x^{\mathrm{T}}.
$$
\end{itemize}
\end{exa}

Example \ref{e1.1.4} shows that the general Definition \ref{d1.1.2} for hypermatrices covers the classical matrix case.

For the case of  general hypermatrices, people are more interested in the case that the eigenvector $x=z^r$, where $z\in \F^n$. We give the following example to depict this.

\begin{exa}\label{e1.1.5} Given a  hypermatrix ${\cal A}=(a_{i_1,i_2,\cdots, i_d})\in \F^{\overbrace{n\times \cdots\times n}^d}$.
${\bf j}=\{j_1,\cdots,j_r\}\subset {\bf i}=\{i_1,\cdots,i_d\}$. We look for an eigenvector of the form $x=z^r$. where $z\in \F^n$, and ${\cal B}:z^r\mapsto z^{d-r}$. Then the eigeneequation becomes
\begin{align}\label{1.1.11}
{\cal A}\times_{{\bf j}} z^r=\lambda z^{d-r}.
\end{align}
\end{exa}

Recently, in addition to the original matrix-matrix STP, the matrix-vector STP is  proposed \cite{che19c} and applied to dimension-varying dynamic (control) systems \cite{che19b}. In addition, the STP of hypermatrices is also proposed \cite{che23}.

\section{Preliminaries}

This section provides some necessary mathematical preliminaries, including (i) STP; (ii) Hypermatrix; (iii) Hypervector .

\subsection{STP of Matrices}

This subsection is a brief review on the STP of matrices \cite{che11,che12}.

\begin{dfn}\label{d2.1.1} Let $A\in {\cal M}_{m\times n}$, $B\in {\cal M}_{p\times q}$, and $t=\lcm(n,p)$. The STP of $A$ and $B$ is defined as follows.
\begin{align}\label{2.1.1}
		A\ltimes B=\left(A\otimes I_{t/n}\right) \left(B\otimes I_{t/p}\right).
	\end{align}
\end{dfn}

The STP is a generalization of the classical matrix product, i.e., when $n=p$,  $A\ltimes B=AB$. Because of this, we  omit the symbol $\ltimes$ in most cases.

One of the most important advantages of STP is that it keeps most of the properties of the classical matrix product. In the following, some basic properties are reviewed, which will be used in the sequel.

\begin{prp}\label{p2.1.2}
	\begin{itemize}
		\item[(i)] (Associativity) Let $A,B,C$ be three matrices of arbitrary dimensions. Then $(A\ltimes B)\ltimes C=A\ltimes (B\ltimes C)$.
		\item[(ii)] (Distributivity) Let $A,B$ be two matrices of same dimension, $C$ is of arbitrary dimension. Then $(A + B)\ltimes C=A\ltimes C + B\ltimes C$, $C\ltimes (A + B)=C\ltimes A + C\ltimes B$.
		\item[(iii)] Let $A,B$ be two matrices of arbitrary dimensions. Then
		\begin{align}\label{2.1.4}
			(A \ltimes B)^T=B^T\ltimes A^T.
		\end{align}
		\item[(iv)] Let $A,B$ be two invertible matrices of arbitrary dimensions. Then $(A \ltimes B)^{-1}=B^{-1}\ltimes A^{-1}$.
	\end{itemize}
\end{prp}

\begin{prp}\label{p2.1.3}
	Let $x\in \F^n$ be a column vector. Then
		\begin{align}\label{2.1.6}
			xA=(I_n\otimes A)x.
		\end{align}
		Let $\omega\in \F^n$ be a row vector. Then $A\omega=\omega(I_n\otimes A)$.
\end{prp}

Note that if $x_i\in \F^i$ $i\in [1,r]$ are column vectors, then by definition,
$\ltimes_{i=1}^rx_i=\otimes_{i=1}^rx_i$.
Then we have the following proposition.

\begin{prp}\label{p2.1.4} Let $T_i\in {\cal M}_{m_i\times n_i}$, $x_i\in \F_{n_i}$, $i\in[1,r]$. Then
	\begin{align*}
		\left(\otimes_{i=1}^rM_i\right)\ltimes_{i=1}^rx_i=\ltimes_{i=1}^r(M_ix_i).
	\end{align*}
\end{prp}

Recall that the contract product of hypermatrices \cite{lim13} is essentially a multi-linear map.  The STP is convenient in expressing such maps. For instance, the product (\ref{1.1.1}), firstly proposed by Qi \cite{qi05}, can be written as
\begin{align}\label{2.1.9}
Mx^{d-1}=M\ltimes \overbrace{x \ltimes \cdots \ltimes x}^{d-1},
\end{align}
where $M$ is understood as the matrix expression of a hypermatrix. Hereafter, $Mx^r$ is always considered as the STP of $M$ with $r$ tuples of $x$.

\subsection{Hypermatrices}

\begin{dfn}\label{d3.1} A hypermatrix of order $d$ and dimension $\{n_1,\cdots,n_d\}$ is a set of data over $\F$.
\begin{align*}
		{\cal A}=\{a_{i_1,\cdots,i_d}\in \F\;|\;i_s\in [1,n_s],\; s\in [1,d]\}.
	\end{align*}
	Denote the set of hypermatrices of order $d$ and dimension $\{n_1,\cdots,n_d\}$ by $\F^{n_1\times \cdots\times n_d}$.

If $n_i=n$, $i\in [1,d]$, ${\cal A}$ is called an equilateral hypermatrix.
\end{dfn}

Let
${\bf  i}=\{i_1,\cdots,i_d\}$,  ${\bf j}=\{j_1,\cdots,j_r\}$, ${\bf k}=\{k_1,\cdots,k_s\}$, and
$$
{\bf i}={\bf j} \bigcup {\bf k}
$$
be a partition, where $r+s=d$. Then ${\cal A}$ can be arranged into a matrix $M^{{\bf k}\times {\bf j}}({\cal A})=(a_{\a,\b})$ of dimension $n_{|{\bf k}|} \times n_{|{\bf j}|}$ with
$$
n_{{\bf j}}=\prod_{i=1}^r n_{j_i},\quad n_{{\bf k}}=\prod_{i=1}^s n_{k_i},
$$
while its columns are indexed by $\{j_1,\cdots,j_r\}$, its rows are indexed by $\{k_1,\cdots,k_s\}$, and the elements are arranged in alphabetic order.

In particular, if ${\bf k}=\emptyset$, then ${\cal A}$ is arranged into a column vector, in the form of
\begin{align}\label{3.1}
V({\cal A})=M^{ {\bf i}\times \emptyset}({\cal A})=(a_{1,\cdots,1}, a_{1,\cdots,2},\cdots, a_{n_1,\cdots,n_d})^T.
\end{align}

When ${\bf j}=\{i_1\}$, it is denoted by
\begin{align}\label{3.2}
	\begin{array}{l}
	M({\cal A})=M^{\{i_1\}\times \{i_2,\cdots,i_d\}}({\cal A})
		=\begin{bmatrix}
			a_{1,\cdots,1}&a_{1,\cdots,2}&\cdots&a_{1,n_2,\cdots,n_d}\\
			a_{2,\cdots,1}&a_{2,\cdots,2}&\cdots&a_{2,n_2,\cdots,n_d}\\
			\vdots&~&~&~\\
			a_{n_1,\cdots,1}&a_{n_1,\cdots,2}&\cdots&a_{n_1,n_2,\cdots,n_d}\\
		\end{bmatrix}
	\end{array}
\end{align}

We give an example to describe this.
\begin{exa}\label{e3.2}
	Given ${\cal A}=[a_{i_1,i_2,i_3}]\in \F^{2\times 3\times 2}$. Then
		$$
	\begin{array}{l}
			V({\cal A})=M^{ \{i_1,i_2,i_3\}\times \emptyset}({\cal A})=[a_{111},a_{112},a_{121},a_{122},
			a_{131},a_{132},a_{211},a_{212},a_{221},a_{222},a_{231},a_{232}]^T.\\
			M({\cal A})=M^{\{i_1\}\times \{i_2,i_3\}}({\cal A})=
			\begin{bmatrix}
				a_{111}&a_{112}&a_{121}&a_{122}&a_{131}&a_{132}\\
				a_{211}&a_{212}&a_{221}&a_{222}&a_{231}&a_{232}
			\end{bmatrix}.\\
		M^{\{i_1,i_3\}\times \{i_2\}}({\cal A})=
	\begin{bmatrix}
			a_{111}&a_{121}&a_{131}\\
			a_{112}&a_{121}&a_{132}\\
			a_{211}&a_{221}&a_{231}\\
			a_{212}&a_{221}&a_{232}\\
		\end{bmatrix}.
	\end{array}
		$$
\end{exa}

Using matrix expression of hypermatrix, the contraction product in the Definition \ref{d1.1.2} for eigenvalue-eigenvector can be expressed as a conventional matrix product. The following proposition shows this. It is an immediate consequence of the definitions and the matrix/vector expressions of the corresponding hypermatrices.

\begin{prp}\label{p3.3}  Let ${\cal A}=(a_{i_1,\cdots,i_d})\in \F^{n_1\times \cdots\times n_d}$,
$$
\begin{array}{l}
{\bf j}=\{j_1,\cdots,j_r\}\subset {\bf i}=\{i_1,\cdots,i_d\},\\
{\bf k}=\{k_1,\cdots,k_s\}={\bf i}\backslash {\bf j},
\end{array}
$$
where $1\leq r\leq d-1$ and $r+s=d$. $x\in \F^{n_{j_1}\times \cdots\times n_{j_r}}$, $y\in \F^{n_{k_1}\times \cdots\times n_{k_s}}$. Denote
$$
A:=M^{{\bf k}\times {\bf j}}({\cal A}).
$$
Assume $y={\cal A}\times_{{\bf j}} x$, then
\begin{align}\label{3.3}
V(y)=AV(x).
\end{align}
\end{prp}

\begin{rem}\label{r3.4} Since any matrix expression of a hypermatrix can uniquely determine the hypermatrix, in the eigenequation (\ref{1.1.10}) and for fixed ${\bf j}\subset {\bf i}$, we can replace ${\cal A}$ by $A$, and use $x$ for $V_x$ and $y={\cal B}(x)$ for $V_y$. Then (\ref{1.1.10}) can be replaced by
\begin{align}\label{3.4}
Ax=\lambda y,
\end{align}
where
\begin{align}\label{3.5}
y={\cal B}(x).
\end{align}
Note that in
(\ref{3.4}) the vectors, matrix, and the product are all classical as defined in linear algebra,
\end{rem}

\begin{rem}\label{r3.5} Recall the existing eigenproblems introduced in Introduction. One sees easily that
\begin{itemize}
\item[(i)] In (\ref{1.1.2})-(\ref{1.1.4}) ${\cal A}\in \R^{\overbrace{n,\cdots, n}^d}$ and
$
M_A=M^{\{i_1\}\times\{i_2,\cdots,i_d\}}$.
\item[(ii)] In (\ref{1.1.5})-(\ref{1.1.6}) ${\cal A}\in \R^{\overbrace{n,\cdots, n}^{2d-2}}$ and
$
M_A=M^{\{i_1,\cdots,i_{d-1}\}\times\{i_d,\cdots,i_{2d-2}\}}$.
\item[(iii)] In (\ref{1.1.7}) ${\cal A}\in \R^{\overbrace{n,\cdots,n}^6}$ and
$
M_A=M^{\{i_1,i_2,i_3\}\times\{i_4,i_5,i_6\}}$.
\end{itemize}
\end{rem}

{\bf Hereafter, we always assume the sub-index ${\bf j}\subset {\bf i}$ is pre-assigned. Then the eigenproblem for a hypermatrix can always be expressed by (\ref{3.4})-(\ref{3.5})}.

\subsection{Hypervectors}

\begin{dfn}\label{d4.1}
\begin{itemize}
\item[(i)] Let $x_i\in \F^{n_i}$, $i\in [1,r]$.
$x=\ltimes_{i=1}^rx_i$ is called a hypervector of degree $r$ and dimension $n_1\times \cdots \times n_r$.
The set of hypervector of degree $r$ and dimension $n_1\times \cdots \times n_r$ is denoted by
$\F^{n_1\ltimes \cdots \ltimes n_r}$.
\item[(ii)] Let $x\in \F^n$, where $n=\prod_{i=1}^rn_i$, if there exist $x^i\in \F^{n_i}$, $i\in [1,r]$ such that
$x=\ltimes_{i=1}^rx_i$, then $x$ is said to be decomposable with respect to $(n_1,\cdots,n_k)$. If there is no confusion,
$x$ is briefly said to be decomposable.

\end{itemize}
\end{dfn}

Let $A\in \F^{n_1\times \cdots \times n_r}$. Then $V(A)$ is a vector. Now consider $x\in \F^{n_1\ltimes \cdots \ltimes n_r}$. Then $x$ can be considered as its vector form $V(x)$ of a hypermatrix $x\in  \F^{n_1\times \cdots \times n_r}$. From this perspective  $x$ can be ``imbedded" into $\F^{n_1\times \cdots \times n_r}$. We conclude that
\begin{align*}
\F^{n_1\ltimes \cdots \ltimes n_r}\subset \F^{n_1\times \cdots \times n_r}.
\end{align*}

When $x=\ltimes_{i=1}^rx_i$, $\{x_i\;|\;i\in [1,r]\}$ are called the components of $x$. An interesting question is: can $x$ uniquely determine its components, or are there multiple ways to express a hypervector by its components? The following proposition provides an answer to it.

\begin{prp}\label{p4.2}
Assume that $x=\ltimes_{i=1}^rx_i=\ltimes_{i=1}^rz_i\neq 0$, $x_i,z_i\in \F^{n_i}$, $i\in [1,r]$.
Then $z_i=k_ix_i,\quad k_i\in \F,\; i\in[1,r]$, and $\prod_{i=1}^rk_i=1$.
\end{prp}

\noindent{\it Proof.} Since $x\neq 0$, all $x_i\neq 0$, $i\in [1,r]$. Assume the $j_i$ th component of $x_i$ neq 0, that is,
$x_i^{j_i}\neq 0$, $i\in [1,r]$.
Set
\begin{align}\label{4.1}
P_i:=\left(\otimes_{k=1}^{i-1}[\d_{n_k}^{j_k}]^T\right)\otimes I_{n_i}\otimes\left([\otimes_{k=i+1}^{r}\d_{n_k}^{j_k}]^T\right),\quad i\in [1,r].
\end{align}
Using Proposition \ref{p2.1.4}, we have
\begin{align*}
P_ix=\left(\prod_{k=1}^{i-1} x_k^{j_k}\prod_{k=i+1}^{r}x_k^{j_k}\right)x_i:=c_ix_i,
\end{align*}
where $c_i\neq 0$. Similarly, $P_iz=d_iz_i$, where $d_i\neq 0$.
Then the conclusion follows.

\hfill $\Box$

\begin{cor}\label{c4.3} Let $x=\ltimes_{i=1}^rx_i\neq 0$, where $x_i\in \F^n$. Then the components of $x$ span a subspace of $\F^n$. That is, $S(x)=\Span\{x_1,\cdots,x_r\}$
is uniquely determined.
\end{cor}

\section{U-Eigenproblem of  Hypermatrices}

\subsection{Monic Decomposition Algorithm}

\begin{dfn}\label{d5.1.1} Let $0\neq x\in \R^n$.
\begin{itemize}
\item[(i)] The Index of $x$, denoted by $\mu(x)$, is defined by
$
\mu(x):=\min\{i\;|\;x_i\neq 0\}.
$
\item[(ii)] The index value of $x$, denoted by $c_0(x)$, is defined by $c_0(x)=x_{e}$, where $e=\mu(x)$.
\item[(iii)] $x$ is called a monic vector if $c_0(x)=1$.
\end{itemize}
\end{dfn}

Note that for $x\neq 0$ $x_0=\frac{1}{c_0(x)}x$ is monic.

Assume $n=\prod_{i=1}^rn_i$. Let $1\leq e\leq n$. By long division, $e$ can be expressed into the following form uniquely.
\begin{align}\label{5.1.1}
e-1=(\cdots(c_1n_2+c_2)n_3+\cdots+c_{r-1})n_r+c_r,
\end{align}
where $0\leq c_j\leq n_{j+1}-1$, $j\in [1,r-1]$.

A straightforward computation shows the following.

\begin{prp}\label{p5.1.2} Consider $x=\ltimes_{i=1}^r x_i$, where $x_i\in \F^{n_i}$, $i\in [1,r]$. Then $x\in \F^n$ with $n=\prod_{i=1}^rn_i$. Assume $\mu(x)=e$, and $e$ is expressed as in (\ref{5.1.1}. Then
\begin{align}\label{5.1.2}
\mu(x_i)=c_i+1,\quad i\in [1,r].
\end{align}
\end{prp}

For our purpose, hereafter we assume  $n_i=n$, $i\in [1,r]$. To calculate the components' indexes from the product index or vise versa, the following provides the formulas,  which are easily verifiable by straightforward computations.

\begin{prp}\label{p5.1.3} Let $x=\ltimes_{1}^rx_i$, where $x_i\in \F^n$, $i\in [1,r]$, $\mu(x)=e$, and $\mu(x_i)=e_i$, $i\in [1,r]$.
\begin{itemize}
\item[(i)] Set $a_0=e-1$, and calculate the following recursively:
\begin{align}\label{5.1.3}
a_s=a_{s+1}n+c_{r-s}, \quad s\in [0,r-1],
\end{align}
where $c_{r-s}\leq n-1$. Then
\begin{align}\label{5.1.4}
\mu(x_i)=c_i+1,\quad i\in[1,s].
\end{align}
\item[(ii)]
\begin{align}\label{5.1.5}
e=\dsum_{i=1}^r(e_i-1)n^{n-i}+1.
\end{align}
\item[(iii)] Assume $x_i=z$, $i\in[1,s]$ and $\mu(z)=e_0$. Then
\begin{align}\label{5.1.6}
(e-1)(n-1)=(e_0-1)(n^r-1).
\end{align}
\end{itemize}
\end{prp}

Define an index-based mapping as
\begin{align}\label{5.1.3}
\Xi^e_{(i;n)}:=[\d_{n_1}^{c_1+1}]^{\mathrm{T}} \otimes \cdots  \otimes [\d_{n_{i-1}}^{c_{i-1}+1}]^{\mathrm{T}}\otimes I_{n_i}\otimes
[\d_{n_{i+1}}^{c_{i+1}+1}]^{\mathrm{T}} \otimes \cdots  \otimes [\d_{n_{r}}^{c_{r}+1}]^{\mathrm{T}},\quad i\in [1,r].
\end{align}

The following Proposition shows how to calculate the component vectors from a hypervector.

\begin{prp}\label{p5.1.3} Assume $x\in \R^n$ with $\mu(x)=e$ and $c_0(x)=a$, and $x$ is decomposable. Then
\begin{align}\label{5.1.4}
\Xi^e_{(i;n)}(x)=ax_i,\quad i\in [1,r].
\end{align}
\end{prp}

Next, we consider how to solve an index-depending equation.
Consider an index-depending linear system
\begin{align}\label{5.1.5}
M_ex=0,
\end{align}
where $e=\mu(x)$. This is a homogeneous equation, so we can, w.l.g., assume $x$ is monic.

To solve $x$ we start from the assumption  $\mu(x)=1$. Then the corresponding $x$ becomes $(1,x_2,\cdots,x_n)^{\mathrm{T}}$.
Then the coefficient matrix becomes a determinant one, say $M_1$. The equation is solvable. Next, we assume $\mu(x)=2$,
Then the corresponding $x$ becomes $(0,1,x_3,\cdots,x_n)^{\mathrm{T}}$. The equation becomes $M_2x=0$. Continuing this process, eventually all the solutions can be obtained.

\begin{rem}\label{r5.1.4}
\begin{itemize}
\item[(i)] Equation (\ref{5.1.4}) is called the Decomposition Formula, which is used to calculate the components of a hypervector.
\item[(ii)] Using Decomposition Formula, one sees easily that if $x$ is monic, then all its component vectors are monic, and vise versa.
\item[(iii)] The method for solving equation (\ref{5.1.5}) will be used to solve eigenequation, where $x=\ltimes_{i=1}^rx_i$ is a hypervector. Then the decomposition formula with the monic algorithm for solving index-depending linear equation is called the {\bf Monic Decomposition Algorithm}.
\item[(iv)] Particularly, if $x$ is diagonal, i.e., $x=z^r$, then $e=\mu(x)$ can only have $n$ possibility (where $z\in \F^n$), which are
\begin{align}\label{5.1.6}
e=(e_0-1)\frac{n^r-1}{n-1}+1,\quad e_0\in [1,n].
\end{align}
Then for solving (\ref{5.1.5}), instead of $n^r$ cases, we have only to consider $n$ cases.
\end{itemize}
\end{rem}

Finally, to ensure our type ${\cal B}$ in eigenproblem being a mapping from subspace to subspace, we need a rigorous description.

\begin{dfn}\label{d5.1.5} Let ${\cal B}:\F^{\overbrace{n\ltimes \cdots\ltimes n}^r}\ra \F^{\overbrace{n\ltimes \cdots\ltimes n}^s}$.
${\cal B}$ is a type-based multi-linear mapping, if
\begin{itemize}
\item[(i)] There exist a set of $s$ matrices $B_j\in {\cal M}_{n\times n^r}$, $j\in [1,s]$, such that for each
 $x=\ltimes_{i=1}^rx_i\in \F^{\overbrace{n\ltimes \cdots\ltimes n}^r}$
\begin{align}\label{5.1.7}
y_j=B_jx,\quad j\in [1,s].
\end{align}
\item[(ii)]
\begin{align}\label{5.1.8}
{\cal B}(x):=\ltimes_{j=1}^sy_j.
\end{align}
\end{itemize}
\end{dfn}

Now the eigenproblem for hypermatrices can be considered as a mapping from $S(x)$ to $S(y)$. By designing type, these two subspaces of $\F^n$ can have required relationship.

\subsection{Type-Based U-Eigenproblem}

\begin{dfn}\label{d5.02.1} Given an equilateral hypermatrix ${\cal A}=(a_{i_1,\cdots,i_d})\in \F^{\overbrace{n\times \cdots \times n}^d}$. Assume the followings:
\begin{itemize}
\item[(i)] There is a preassigned index partition
$$
{\bf i}={\bf j}\bigcup {\bf k}
$$
with $|{\bf j}|=r$, $1\leq r\leq d-1$, $|{\bf k}|=s$, $r+s=d$.

\item[(ii)] There is a monic
$$
x=\ltimes_{j=1}^rx_j\in \F^{\overbrace{n\ltimes \cdots \ltimes n}^r},
$$
where the entries of $x$ are labeled by ${\bf j}=(j_1,\cdots,j_r)$ as
\begin{align}\label{5.02.1}
x=\left\{\prod_{i=1}^r x_i^{j_i}:=x_{j_1,j_2,\cdots,j_r}\;|\; j_i\in [1,n], i\in [1,r]\right\},
\end{align}
and $x_{i}^{j_i}$ is the $j_i$-th entry of $x_{i}$.

\item[(iii)] There is a type-based multi-linear mapping ${\cal B}\in \F^{\overbrace{n\ltimes \cdots\ltimes n}^r}\ra \F^{\overbrace{n\ltimes \cdots\ltimes n}^s}$, and a scale $\lambda$, such that
\begin{align}\label{5.02.2}
{\cal A}\times_{{\bf j}}x=\lambda {\cal B}x.
\end{align}
\end{itemize}
Then $x$ is called a U-eigenvector of ${\cal A}$ w.r.t. type ${\cal B}$ and eigenvalue $\lambda$.
\end{dfn}

\begin{dfn}\label{d5.02.2}
\begin{itemize}
\item[(i)] $x=\ltimes_{i=1}^rx_i$ is said to be a diagonal hypervector, if $x_i=z$, $i\in [1,r]$.
\item[(ii)] Assume $x$ is a U-eigenvector (of ${\cal A}$ w.r.t. type ${\cal B}$ and eigenvalue $\lambda$).
In addition, both $x$ and $y$ are diagonal, $x$ is a D-eigenvector.
\end{itemize}
\end{dfn}

\begin{rem}\label{r5.02.3}
\begin{itemize}
\item[(i)] U-eigenvector has significant geometric meaning.  Both $x$ and $y$ represent the corresponding subspaces $S(x)$ and $S(y)$. Type ${\cal B}$ describing the relationship between $S(x)$ and its image $S(y)$. Say, $S(y)\subset S(x)$, or $S(y)\perp S(x)$, (both will be described in the sequel).
\item[(ii)] Based on the subspace perspective, $x$ and $kx$ represent the same subspace. Unfortunately,
assume $x$ satisfies (\ref{5.02.2}), $kx$ does not. According to (\ref{5.02.2}) it is easy to see that
$$
{\cal A}(kx)=k{\cal A}(x),\quad
{\cal B}(kx)=k^s{\cal B}(x),
$$
which leads to the conclusion that $kx$ is the H eigenvector w.r.t. eigenvalue $\lambda k^{s-1}$. Though, $x$ and $kx$ represent the same eigen-subspace, they could correspond different eigenvalues. To make the eigenvalue unique, we demand the {\bf eigenvector $x$ to be monic}.

\item[(iii)] From (ii) one sees that for hypermatrix eigenvector is more fundamental than eigenvalue. $x$ and $kx$ are both eigenvectors, which represent same subspace, both they correspond to different eigenvalues.
\end{itemize}
\end{rem}

The following is a fundamental result, which converts the U-eigenproblem to a linear  algebraic equation.

\begin{thm}\label{t5.02.4} The U-eigenproblem can be solved by solving a general matrix eigenequation.
\end{thm}

\noindent{\it Proof.} Set
$$
A=M^{{\bf k}\times {\bf j}}({\cal A}),
$$
Then we have
$$
{\cal A}\times_{{\bf j}}x=Ax,
$$
where $x=\ltimes_{i=1}^rx_i$.

Now consider ${\cal B}(x)$. By definition
$$
{\cal B}(x)=\ltimes_{j=1}^sy^j=\ltimes_{j=1}^s(B_ix).
$$
Using Proposition \ref{p2.1.3}, we have
$$
\begin{array}{l}
\ltimes_{j=1}^s(B_ix)=B_1(I_{n^r}\otimes B_2)x^2(B_3x)\cdots(B_sx)\\
=\cdots=B_1(I_{n^r}\otimes B_2)\cdots (I_{n^{(s-1)r}}\otimes B_s)x^s\\
:=\tilde{B}x^s.
\end{array}
$$
Then (\ref{5.02.2}) becomes
$$
Ax=\lambda \tilde{B}x^s.
$$
Expressing $x=\Xi^e_{(1, n^r)}x^s$ and setting $\tilde{A}=A\Xi^e_{(1, n^r)}$, we have
\begin{align}\label{5.02.3}
\tilde{A}\xi =\lambda \tilde{B}\xi,
\end{align}
where $\xi=x^s$.

Now (\ref{5.02.3}) is a general matrix eigenequation.  Solving it yields candidate $\xi$. Then Monic Decomposition Algorithm can be used to pick out decomposable $\xi$, which yield $x$. Using Monic Decomposition Algorithm again, the solution $\{x_1,\cdots,x_r\}$ can finally be obtained.
\hfill $\Box$

Solving general matrix eigenequation will be discussed in next section. The numerical solution can be obtained via MatLab.

\subsubsection{Eigenequation for D-Eigenproblem}

Since D-eigenproblem is relatively easy and it is the object of all literature so far, we discuss it first.

 For D-eigenproblem same procedure leads to the  algebraic expression of (\ref{5.2.2}) as
\begin{align}\label{5.2.4}
Az^r=\lambda z^s.
\end{align}
Making it homogeneous, we
define
$$
E=\begin{cases}
I_{n^s}\otimes \left[(\d_n^{\mu(x)})^{\mathrm{T}}\right]^{r-s},\quad r>s,\\
I_{n^r}\otimes \left[(\d_n^{\mu(x)})^{\mathrm{T}}\right]^{s-r},\quad r<s.\\
\end{cases}
$$
A straightforward computation shows the following.

\begin{cor}\label{c5.2.5}  The linear algebraic equation for D-eigenproblem is
\begin{align}\label{5.2.5}
\begin{cases}
\tilde{A}z^s=(AE)z^s=\lambda z^s,\quad r<s,\\
Az^r=\lambda Ez^r.
\end{cases}
\end{align}
\end{cor}

Next example shows the eigenproblems mentioned in Section 1 are all D-eigenproblem.

\begin{exa}\label{e5.2.6} Consider

\begin{enumerate}
		\item[] Case 1: $r=s$ (then $d$ is even.)
	\end{enumerate}
	
\begin{itemize}
\item[(i)] The type  mapping ${\cal B}$ defined by  (\ref{1.1.2}) is of this style.
Set
$$
d:=n^{r-1}+n^{r-2}+\cdots +n,
$$
and
$$
B^{\mathrm{T}}=\d_{n^{r}}[1,d+2,2d+3,\cdots,(n-1)d+n].
$$
Then
$$
Bz^r=(z_1^r,z_2^r,\cdots, z_n^r)^T.
$$
To be specific, we assume $r=d-1=3$, $n=2$ and $s=1$, then we have
\begin{align}\label{5.2.1}
x=z^3=\begin{bmatrix}
z_1^3\\
z_1^2z_2\\
z_1z_2z_1\\
z_1z^2_2\\
z_2z_1^2\\
z_2z_1z_2\\
z_2^2z_1\\
z_2^3
\end{bmatrix}
\end{align}
Then
$$
B^{\mathrm{T}}=\d_8[1,8].
$$

\item[(ii)]  The type  mapping ${\cal B}$ defined by  (\ref{1.1.5}) is of this style.
To be specific, we assume $r=d-1=3$, $n=2$ and $s=1$, then
$$
(z_1+z_2)^2\begin{bmatrix}
z_1\\
z_2
\end{bmatrix}=
\begin{bmatrix}
z_1^3+2z_1^2z_2+z_1z_2^2\\
z_1^2z_2+2z_1z_2^2+z_2^3
\end{bmatrix}.
$$
Using (\ref{5.2.1}), we have
$$
B=\begin{bmatrix}
1&2&0&0&0&0&1&0\\
0&1&0&2&0&0&0&1\\
\end{bmatrix}.
$$

\item[(iii)]  The type  mapping ${\cal B}$ defined by  (\ref{1.1.6}) is of this style.
The matrix $B$ has already been included in (\ref{1.1.6}).

\item[(iv)]  The type  mapping ${\cal B}$ defined by  (\ref{1.1.7}) is of this style.

We have $r=d-1=3$, $s=1$, and assume $n=2$, then
$$
(z_1,z_2)\begin{bmatrix}
z_1\\
z_2
\end{bmatrix}z=
\begin{bmatrix}
z_1^3+z_1z_2^2\\
z_1^2z_2+z_2^3
\end{bmatrix}.
$$
Using (\ref{5.2.1}), we have
$$
B=\begin{bmatrix}
1&0&0&1&0&0&0&0\\
0&1&0&0&0&0&0&1\\
\end{bmatrix}.
$$
\end{itemize}
\end{exa}

\begin{itemize}
\item[] Case 2:
\end{itemize}

Assume $x=z^r$ and ${\cal B}(x)=Bz^s$, where $s\leq r$. (Particularly, $s=1$.)

\begin{exa}\label{e5.2.2}
The following facts are obvious.
\begin{itemize}
\item[(i)] The type mapping ${\cal B}$ defined by  (\ref{1.1.3}) is of this style, (where $s=1$).
\item[(ii)] The type mapping ${\cal B}$ defined by  (\ref{1.1.4}) is of this style, (where $s=1$) .
\end{itemize}
\end{exa}

\begin{itemize}
\item[] Case 3:
\end{itemize}

Assume $x=z^r$ and ${\cal B}(x)=Bz^s$, where $s > r$.

The following example shows this case.

\begin{exa} \label{e5.2.3} \cite{qi05}
\begin{itemize}
\item[(i)] Assume ${\cal A}\in \F^{\overbrace{n\times \cdots\times n}^{2m-1}}$ and $A=M^{{\bf k}\times {\bf j}}(A)$, where ${\bf i}=\{i_1,\cdots,i_{2m-1}\}$,
${\bf j}=\{i_1,\cdots,i_m\}$, ${\bf k}=\{i_{m+1},\cdots,i_{2m-1}$. The eigenproblem is defined as follows.
\begin{align}\label{5.2.2}
Ax^{m-1}=\lambda [x^{\mathrm{T}}x]^{m-2}x.
\end{align}

\item[(ii)] A numerical example in \cite{qi07} is as follows:

\begin{align}\label{5.2.3}
\begin{cases}
x_1^2=\lambda x_0x_1,\\
x_2^2=\lambda x_0x_2,\\
x_0^2=x_1^2-x_2^2.
\end{cases}
\end{align}
Set $x=(x_1,x_2)^{\mathrm{T}}$, then a straightforward computation converts (\ref{5.2.3}) into an eigenequation as
\begin{align}\label{5.2.4}
Ax^2=\lambda Bx^3,
\end{align}
where
$$
A=\begin{bmatrix}
1&0&0&0\\
0&0&0&1
\end{bmatrix},\quad
B=\begin{bmatrix}
1&0&0&1&0&0&0&0\\
0&1&0&0&0&0&0&1
\end{bmatrix}.
$$
\end{itemize}
\end{exa}

\subsubsection{Eigenequation for U-Eigenproblem}

The constructive proof of Theorem \ref{t5.02.4} provides a process to construct the eigenequation for U-eigenproblem. We give an example to describe it.

\begin{exa}\label{e5.2.4}
Given ${\cal A}=(a_{i_1,i_2,i_3,i_4})\in \F^{2\times 2\times 2 \times 2}$, where
$$
a_{i_1,i_2,i_3,i_4}=
\begin{cases}
1,\quad {\bf i}=\{(1,1,1,1), (1,1,2,1), (2,1,2,2),(1,2,1,2),(2,2,2,1)\},\\
2,\quad {\bf i}=\{(1,1,1,2),(2,1,1,2),\}\\
0,\quad \mbox{Otherwise}.
\end{cases}
$$
Then
$$
\begin{array}{l}
A=M^{(i_1,i_2)\times (i_3,i_4)}({\cal A})\\
=\begin{bmatrix}
1&2&1&0\\
0&0&0&1\\
0&2&0&1\\
0&0&1&0\\
\end{bmatrix}.
\end{array}
$$
Assume
$$
\begin{array}{l}
x=x_1x_2,\\
y_1=x_1+x_2,\\
y_2=x_1,\\
y=y_1y_2.
\end{array}
$$
Then
$$
\begin{array}{l}
y_1=(I_2\otimes [\d_2^{e_2}]^{\mathrm{T}})x+( [\d_2^{e_1}]^{\mathrm{T}}\otimes I_2)x:=(B_1+B_2)x,\\
y_2=(I_2\otimes [\d_2^{e_2}]^{\mathrm{T}})x:=B_1x.
\end{array}
$$
\begin{itemize}
\item[(i)] $e_1=e_2=1$:
$$B_1=I_2\otimes  [\d_2^{1}]^{\mathrm{T}}=\begin{bmatrix} 1&0&0&0\\0&0&1&0\end{bmatrix}.
$$
$$B_2= [\d_2^{1}]^{\mathrm{T}}\otimes I_2=\begin{bmatrix} 1&0&0&0\\0&1&0&0\end{bmatrix}.
$$
\item[(ii)] $e_1=1,~e_2=2$:
$$B_1=I_2\otimes  [\d_2^{2}]^{\mathrm{T}}=\begin{bmatrix} 0&1&0&0\\0&0&0&1\end{bmatrix}.
$$
$$B_2= [\d_2^{1}]^{\mathrm{T}}\otimes I_2=\begin{bmatrix} 1&0&0&0\\0&1&0&0\end{bmatrix}.
$$
\item[(iii)] $e_1=2,~e_2=1$:
$$B_1=I_2\otimes  [\d_2^{1}]^{\mathrm{T}}=\begin{bmatrix} 1&0&0&0\\0&0&1&0\end{bmatrix}.
$$
$$B_2= [\d_2^{2}]^{\mathrm{T}}\otimes I_2=\begin{bmatrix} 0&0&1&0\\0&0&0&1\end{bmatrix}.
$$
\item[(iv)] $e_1=e_2=2$:
$$B_1=I_2\otimes  [\d_2^{2}]^{\mathrm{T}}=\begin{bmatrix} 0&1&0&0\\0&0&0&1\end{bmatrix}.
$$
$$B_2= [\d_2^{2}]^{\mathrm{T}}\otimes I_2=\begin{bmatrix} 0&0&1&0\\0&0&0&1\end{bmatrix}.
$$
\end{itemize}

It follows that
$$
Ax=[B_1(I_4\otimes B_2)]x^2:=Bx^2,
$$

\begin{itemize}
\item[(i)] $e_1=e_2=1$:
$$
\begin{array}{l}
B=\left[
\begin{array}{cccccccccccccccc}
1&0&0&0&0&0&0&0&0&0&0&0&0&0&0&0\\
0&1&0&0&0&0&0&0&0&0&0&0&0&0&0&0\\
0&0&0&0&0&0&0&0&1&0&0&0&0&0&0&0\\
0&0&0&0&0&0&0&0&0&1&0&0&0&0&0&0\\
\end{array}
\right].
\end{array}
$$
\item[(ii)] $e_1=1,~e_2=2$:
$$
\begin{array}{l}
B=\left[
\begin{array}{cccccccccccccccc}
0&0&0&0&1&0&0&0&0&0&0&0&0&0&0&0\\
0&0&0&0&0&1&0&0&0&0&0&0&0&0&0&0\\
0&0&0&0&0&0&0&0&0&0&0&0&1&0&0&0\\
0&0&0&0&0&0&0&0&0&0&0&0&0&1&0&0\\
\end{array}
\right].
\end{array}
$$
\item[(iii)] $e_1=2,~e_2=1$:
$$
\begin{array}{l}
B=\left[
\begin{array}{cccccccccccccccc}
0&0&1&0&0&0&0&0&0&0&0&0&0&0&0&0\\
0&0&0&1&0&0&0&0&0&0&0&0&0&0&0&0\\
0&0&0&0&0&0&0&0&0&0&1&0&0&0&0&0\\
0&0&0&0&0&0&0&0&0&0&0&1&0&0&0&0\\
\end{array}
\right].
\end{array}
$$
\item[(iv)] $e_1=e_2=2$:
$$
\begin{array}{l}
B=\left[
\begin{array}{cccccccccccccccc}
0&0&0&0&0&0&1&0&0&0&0&0&0&0&0&0\\
0&0&0&0&0&0&0&1&0&0&0&0&0&0&0&0\\
0&0&0&0&0&0&0&0&0&0&0&0&0&0&1&0\\
0&0&0&0&0&0&0&0&0&0&0&0&0&0&0&1\\
\end{array}
\right].
\end{array}
$$
\end{itemize}

Then the algebraic equation is expressed as
$$
Ax=\lambda Bx^2.
$$
Note that
$$
Ax=A[I_4\otimes (\d_4^{e_x})^{\mathrm{T}}]x^2:=\tilde{A}x^2.
$$
\begin{itemize}
\item[(i)] $e_1=e_2=1$:
$$
\begin{array}{l}
\tilde{A}=A(I_4\otimes [\d_4^1]^{\mathrm{T}})\\
=\left[
\begin{array}{cccccccccccccccc}
1&0&0&0&2&0&0&0&1&0&0&0&0&0&0&0\\
0&0&0&0&0&0&0&0&0&0&0&0&0&0&0&0\\
0&0&0&0&2&0&0&0&0&0&0&0&1&0&0&0\\
0&0&0&0&0&0&0&0&0&0&0&0&0&0&0&0\\
\end{array}
\right].
\end{array}
$$
\item[(ii)] $e_1=1, ~e_2=2$:
$$
\begin{array}{l}
\tilde{A}=A(I_4\otimes [\d_4^2]^{\mathrm{T}})\\
=\left[
\begin{array}{cccccccccccccccc}
0&1&0&0&0&2&0&0&0&1&0&0&0&0&0&0\\
0&0&0&0&0&0&0&0&0&0&0&0&0&0&0&0\\
0&0&0&0&0&2&0&0&0&0&0&0&0&1&0&0\\
0&0&0&0&0&0&0&0&0&0&0&0&0&0&0&0\\
\end{array}
\right].
\end{array}
$$
\item[(iii)] $e_1=2,~e_2=1$:
$$
\begin{array}{l}
\tilde{A}=A(I_4\otimes [\d_4^3]^{\mathrm{T}})\\
=\left[
\begin{array}{cccccccccccccccc}
0&0&1&0&0&0&2&0&0&0&1&0&0&0&0&0\\
0&0&0&0&0&0&0&0&0&0&0&0&0&0&0&0\\
0&0&0&0&0&0&2&0&0&0&0&0&0&0&1&0\\
0&0&0&0&0&0&0&0&0&0&0&0&0&0&0&0\\
\end{array}
\right].
\end{array}
$$
\item[(iv)] $e_1=e_2=2$:
$$
\begin{array}{l}
\tilde{A}=A(I_4\otimes [\d_4^4]^{\mathrm{T}})\\
=\left[
\begin{array}{cccccccccccccccc}
0&0&0&1&0&0&0&2&0&0&0&1&0&0&0&0\\
0&0&0&0&0&0&0&0&0&0&0&0&0&0&0&0\\
0&0&0&0&0&0&0&2&0&0&0&0&0&0&0&1\\
0&0&0&0&0&0&0&0&0&0&0&0&0&0&0&0\\
\end{array}
\right].
\end{array}
$$

\end{itemize}

We finally have the eigenequation as
\begin{align}\label{5.2.5}
(\tilde{A}-\lambda B)x^2.
\end{align}
\end{exa}

\section{Solving Eigenequation}

From previous section one sees that the eigenequation of U-eigenvalue-eigenvector can be converted into the following form
\begin{align}\label{6.0.1}
(A-\lambda B)x^r=0.
\end{align}
The general solving algorithm we proposed is as follows:
\begin{alg}\label{6.0.2}
\begin{itemize}
\item[] Step 1. Set $\xi=x^r$ and find $\lambda$ such that $K(A-\lambda B)\neq 0$, where $K$ is the kernel set.
\item[] Step 2. For each feasible $\lambda$ find the corresponding solution set $K$ with a basis $\{\xi_1,\cdots, \xi_s\}$.
\item[] Step 3. Verifying if $\xi=\dsum_{i=1}^{s}a_i\xi_i$ is decomposable?
Each decompasable $\xi$ is  a U-eigenvector. Moreover, if $\xi$ is diagonal, it is a D-eigenvector.
\end{itemize}
\end{alg}

\subsection{Solving Essential Eigenvalue}.

\begin{dfn}\label{d6.1.1} Consider
$$
E(\lambda)=(A-\lambda B)\in {\cal M}_{m\times n}.
$$
Define the generic rank of $E(\lambda)$ by
\begin{align}\label{6.1.1}
r_g=\max\{r\;\rank(E(\lambda))\}.
\end{align}
\begin{itemize}
\item[(i)] $\lambda$ is called a $B$-type eigenvalue of $A$, if $\rank (E(\lambda))<n$.
\item[(ii)] If $\lambda$ is an eigenvalue of $A$, and $\rank (E(\lambda))=r_g$, $\lambda$ is called a quasi eigenvalue.
Otherwise,   $\lambda$ is an essential  eigenvalue.
\end{itemize}
\end{dfn}

Since $r_g\leq n$, we can choose $r_g$ rows from $(A-\lambda B$ which is of generic rank $r_g$.
These rows form $\tilde(E)(\lambda)=(\tilde{A}-\lambda \tilde{B})$. Then it is easy to see that $\lambda$ is an essential eigenvalue of $E$, if and only if,  $\lambda$ is an essential eigenvalue of $\tilde{E}$.

Assume $\tilde{B}$ has full row rank, i.e., $\rank(\tilde{B})=r_g$.\footnote{By choosing proper type ${\cal B}$, we can w.l.g. assume its corresponding $B$ is of full rank. That is, $\rank(B)=\min(m,n)$. Then it is easy to find $\tilde{B}$, which satisfies this assumption.}

In the following we (w.l.g) assume the followings:

\noindent{\bf A1}: The generic rank of $A-\lambda B\in {\cal M}_{m\times n}=m$.
And
\noindent {\bf A2}:  $B$ has full row rank, i.e.,  $\rank(B)=m$.

The follows result comes from linear algebra.

\begin{prp}\label{p6.1.2} $\lambda$ is an essential eigenvalue of $A-\lambda B$, if and only if,
\begin{align}\label{6.1.2}
\det\left[(A-\lambda B)(A-\overline{\lambda} B)^{\mathrm{T}}\right]=0.
\end{align}
\end{prp}

\begin{rem}\label{r6.1.3}
\begin{itemize}
\item[(i)] Particularly, when $A$, $B$ are square matrices, (\ref{6.1.2}) becomes
\begin{align}\label{6.1.2}
\det(A-\lambda B)=0.
\end{align}
When $B=I$, we have conventional eigenequation for square matrix. Hence, conventional eigenvalues are essential. Moreover, when $B\neq I$, we have (generalized) $B$-type eigenvalue-eigenvector for square matrices.
\item[(ii)] Set
\begin{align}\label{6.1.3}
\Psi=B^{\mathrm{T}}(BB^{\mathrm{T}})^{-1}.
\end{align}
Then we consider $A\Psi$. It is easy to verify that if $x$ is an eigenvector of $A\Psi$ w.r.t. $\lambda$, then $\Psi x$ is an eigenvector of $A$ w.r.t. type $B$ and eigenvalue $\lambda$.  (mightbe quasi). If $B$ is square (and invertible), then $Psi=B^{-1}$. Moreover, $x$ is an eigenvector of $AB^{-1}$ w.r.t. $\lambda$, if and only if, $Bx$ is a (generalized) eigenvector of $A$ w.r.t. $\lambda$ and type $B$.
\item[(iii)] Given eigenequation $A-\lambda B$. Denote by $q(A,B)$ and $e(A,B)$ the set of quasi and essential eigenvalues of $(A,B)$ respectively. Then $q(A,B)\bigcup e(A,B)=\C$. Moreover, $e(A,B)$ is finite.
\end{itemize}
\end{rem}

We give a numerical example.

\begin{exa}\label{e6.1.4} Given
$$
A=\begin{bmatrix}
1&0&0\\0&1&0\\
\end{bmatrix};\quad
B=\begin{bmatrix}
1&0&0\\0&0&1\\
\end{bmatrix}.
$$
Consider the eigenequation
\begin{align}\label{6.1.4}
A-\lambda B=0.
\end{align}
It is easy to see that $r_g(A-\lambda B)=2$.
\begin{itemize}
\item[(i)] Eigenvalue:
Consider the Kernel space of $A-\lambda B$. It is easy to find that
$$
K(\lambda)=\Span\{0,\lambda, 1\}.
$$
Then any $\lambda\in \C$ is an eigenvalue of $A$ with type $B$. The corresponding eigenvector is: $x=(0,\lambda,1)^{\mathrm{T}}$.

\item[(ii)] Essential Eigenvalue:
Consider equation (\ref{6.1.2}). Let $\lambda =a+bi$. A straightforward computation simplifies (\ref{6.1.2}) to
$$
(1-2a+a^2+b^2)(1+a^2+b^2)=0.
$$
The solution is: $a=1$ and $b=0$. Hence, the only essential eigenvalue is $\lambda=1$.

\item[(iii)] Set $\Psi=B^{\mathrm{T}}(BB^{\mathrm{T}})^{-1}$ and calculate
$$
A\Psi=\begin{bmatrix}
1&0\\0&0
\end{bmatrix}.
$$
Then it has two eigenvalue: $0$ and $1$, where $0$ is quasi eigenvalue and $1$ is essential eigenvalue.
\end{itemize}
\end{exa}

\subsection{Solving Eigenequation}

To solve eigenequation, we need only to consider essential eigenvalue. In our framework, it becomes a linear algebraic problem.

\begin{exa}\label{e6.2.1} Recall Example \ref{e6.1.4}. The essential eigenvalue is $\lambda=1$. Then
$$
A-\lambda B=
\begin{bmatrix}
0&0&0\\
0&1&-1\\
\end{bmatrix}
$$
Then the solution set is
$$
K_e=\Span\left\{\begin{bmatrix}1\\0\\0\end{bmatrix}, \begin{bmatrix}0\\1\\-1\end{bmatrix}\right\}.
$$
\end{exa}

\begin{rem}\label{r6.2.2} From previous example one sees easily that  why the kernel space of pseudo eigenvalues need not to be considered.
Because the kernel space $K_q$ of quasi eigenvalues is a subspace of the kernel space of essential eigenvalues. So if there is some solutions from $K_q$, it can also be found in any $K_e$. This is extremely important. Because it makes the eigenproblem finitely solvable.
\end{rem}

\subsection{Verifying Decomposibility}

 We use Proposition \ref{p5.1.3} to verify the decomposibility. We can assume $x$ is monic, because  $x$ is decomposible, if and only if, $x=c_0(x)x_0$ and the monic $x_0$ is decomposable. Moreover, $x=z^r$ is the solution of eigenequation, if and only if, $x_0$ is the solution.  First, we use an example to depict the algorithm.

 \begin{exa}\label{e6.3.1} \cite{chi13,li14}
Assume the eigenequation is of the form in (\ref{1.1.5}). Particularly, consider $x=(x_1,x_2)^{\mathrm{T}}$, $d=6$, and $r=s=3$. The equation is expressed as
\begin{align}\label{6.3.1}
\begin{array}{l}
x_1^3=\lambda (x_1+x_2)^2x_1,\\
x_2^3=\lambda (x_1+x_2)^2x_2.\\
\end{array}
\end{align}
Clearly, this is D-eigenproblem.

Then it is easy to convert it into the eigenequation form as
\begin{align}\label{6.3.2}
(A-\lambda B)x^3=0,
\end{align}
where
$$
A=\begin{bmatrix}
1&0&0&0&0&0&0&0\\
0&0&0&0&0&0&0&1\\
\end{bmatrix};\quad
B=\begin{bmatrix}
1&2&0&0&0&0&1&0\\
0&1&0&2&0&0&0&1\\
\end{bmatrix}.
$$
Then
$$
A-\lambda B= \begin{bmatrix}
1-\lambda &-2\lambda&0&0&0&0&-\lambda&0\\
0&-\lambda&0&-2\lambda&0&0&0&1-\lambda\\
\end{bmatrix}.
$$

First, it is easy to calculate the quasi kernel
$$
K_q=\Span\left\{\xi_1,\xi_2,\xi_3,\xi_4,\xi_5,\xi_6\right\},
$$
where
$$
\begin{array}{l}
\xi_1=(0,0,1,0,0,0,0,0)^{\mathrm{T}},\\
\xi_2=(0,0,0,0,1,0,0,0)^{\mathrm{T}},\\
\xi_3=(0,0,0,0,0,1,0,0)^{\mathrm{T}},\\
\xi_4=(0,2,0,-1,0,0,-4,0)^{\mathrm{T}},\\
\xi_5=(4\lambda,2(1-\lambda),0,\lambda-1,0,0,0,0)^{\mathrm{T}},\\
\xi_6=(2\lambda,1-\lambda,0,0,0,0,0,\lambda)^{\mathrm{T}}.\\
\end{array}
$$

Next, we calculate essential eigenvalue. Set $\lambda=a+bi$, a careful calculation shows
$$
\begin{array}{l}
\det\left[(A-\lambda B)\overline{(A-\lambda B)^{\mathrm{T}}}\right]\\
=\left((1-a)^2+a^2+2b^2)\right)\left((1-a)^2+9a^2+10b^2)\right).
\end{array}
$$
So there is no essential eigenvalue.

Now we are ready to verify the dicomposability.

\begin{itemize}
\item[] Case 1: $\lambda=0$
\end{itemize}

Then the vector in $K_q$ can be expressed as

\begin{align}\label{6.3.4}
\xi=\dsum_{i=1}^6 c_i\xi_i
=(0,2c_4+2c_5+c_6, c_1, -c_4-c_5, c_2, c_3, -4c_4,0)^{\mathrm{T}}.
\end{align}

Assume the diagonal solution $x^3$ satisfies (type) $\mu(x)=1$, then $\mu(x^3)=1$, so we need to assume $\mu(x)=2$, then $\mu(x^3)=8$. Comparing with $\xi$ in (\ref{6.3.4}), since $\xi(8)=0$, $\mu(x)=2$ is also not allowed. We conclude that there is no D-eigenvector.

\begin{itemize}
\item[] Case 2: $\lambda\neq 0$
\end{itemize}

We have $\xi\in K$ as
$$
\xi=\dsum_{i=1}^6c_i\xi_i
=\begin{bmatrix}
c_1(4c_5+2c_6)\lambda\\
2c_4+2c_5(1-\lambda) +c_6(1-\lambda)\\
c_1\\
-c_4+c_5(\lambda-1)\\
c_2\\
c_3\\
-4c_4\\
c_6\lambda
\end{bmatrix}
$$

\begin{itemize}
\item[] Case 2.1 $\mu(x)=1$:
\end{itemize}
Denoting $\Xi_i:=\Xi^e_{(i,2)}$, $i=1,2,3$, we have
$$
\begin{array}{l}
\Xi_1=I_2\otimes [\d_4^1]^{\mathrm{T}}
=
\begin{bmatrix}
1&0&0&0&0&0&0&0\\
0&0&0&0&1&0&0&0\\
\end{bmatrix},\\
\Xi_2= [\d_2^1]^{\mathrm{T}}\otimes I_2\otimes [\d_2^1]^{\mathrm{T}}
=
\begin{bmatrix}
1&0&0&0&0&0&0&0\\
0&0&1&0&0&0&0&0\\
\end{bmatrix},\\
\Xi_3= [\d_4^1]^{\mathrm{T}}\otimes I_2
=
\begin{bmatrix}
1&0&0&0&0&0&0&0\\
0&1&0&0&0&0&0&0\\
\end{bmatrix},\\
\end{array}
$$
Then we have
$$
\begin{array}{l}
x_1=\Xi_1\xi
=
\begin{bmatrix}
(4c_5+2c_6)\lambda\\
c_2\\
\end{bmatrix},\\
x_2=\Xi_2\xi
=
\begin{bmatrix}
(4c_5+2c_6)\lambda\\
c_1\\
\end{bmatrix},\\
x_3=\Xi_3\xi
=
\begin{bmatrix}
(4c_5+2c_6)\lambda\\
2c_4+2c_5(1-\lambda)+c_6(1-\lambda)\\
\end{bmatrix}.
\end{array}
$$

To find $D$ eigenvector we set $x_1=x_2=x_3$.
Denote $\theta=c_1=c_2$, and set
$$
(4c_5+2c_6)\lambda =1.
$$
$$
2c_4+2c_5(1-\lambda)+c_6(1-\lambda)=\theta.
$$
Then we have
$$
\begin{bmatrix}
1\\
\theta
\end{bmatrix}^3=\xi.
$$
That is,

\begin{align}\label{6.3.7}
\left\{
\begin{array}{llr}
(4c_5+2c_6)\lambda=1&~&(.1)\\
2c_4+2c_5(1-\lambda) +c_6(1-\lambda)=\theta&~&(.2)\\
c_1=\theta&~&(.3)\\
-c_4+c_5(\lambda-1)=\theta^2&~&(.4)\\
c_2=\theta&~&(.5)\\
c_3=\theta^2&~&(.6)\\
-4c_4=\theta^2&~&(.7)\\
c_6\lambda=\theta^3.&~&(.8)\\
\end{array}\right.
\end{align}

Assume $\theta=0$, it is a feasible case because
$$
\theta=0;~~\mbox{and}~~c_1=c_2=c_3=c_4=c_6=0;\quad \lambda=1;c_5=\frac{1}{4}
$$
is a solution of (\ref{6.3.7}).

Then we know $x=[(1,0)^{\mathrm{T}}]^3$ is a monic D-eigenvector  of ${\cal A}$ w.r.t. eigenvalue $r=1$.

Next, we assume $\theta\neq 0$. Then we have
$$
c_1=c_2=\theta;\quad c_3=\theta^2;\quad c_4=-\frac{1}{4}\theta^2
$$
Form (.4), (.2), and (.8) we have
$$
c_6=c_6\lambda +2\theta^2=\theta^3+2\theta^2,
$$
and
$$
\lambda =\frac{\theta^3}{\theta^3+2\theta^2}=\frac{\theta}{\theta+2}.
$$
Using (.2) again yields
$$
c_5=\frac{-2c_4-c_6+\theta^3}{2(1-\lambda)}=-\frac{3}{8}\theta^2(\theta+2).
$$
Putting the solution into (.1) yields $\theta^2=0$, into (.2) yields $\theta=0$. This contradiction shows that there is no other solution for $\mu(x)=1$.

\begin{itemize}
\item[] Case 2.2 $\mu(x)=2$:
\end{itemize}

In this case the only monic solution is $x=(0,1)^{\mathrm{T}}$. Then we can simply choose
$$
c_6=\lambda=1;\quad c_1=c_2=c_3=c_4=0;\quad c_5=-0.5.
$$
It follows that $x=(0,1)^{\mathrm{T}}$ is a monic D-eigenvector of ${\cal A}$ w.r.t. eigenvalue $\lambda=1$ and type ${\cal B}$.
\end{exa}

 \begin{exa}\label{e6.3.2}

To consider the corresponding U-eigenproblem, we can modify (\ref{6.3.1}) to the following one. Consider
\begin{align}\label{6.3.10}
Axyz=\lambda Bxyz,
\end{align}
where the $A$ and $B$ are the same as in Example \ref{6.3.1}. Then (\ref{6.3.10}) can be expressed into component-wise form as
\begin{align}\label{6.3.11}
\begin{array}{l}
x_1y_1z_1=\lambda(x_1y_1z_1+2x_1y_1z_2+x_2y_2z_1),\\
x_2y_2z_2=\lambda(x_1y_1z_2+2x_1y_2z_2+x_2y_2z_2).\\
\end{array}
\end{align}

Now we can consider 8 different cases: $\mu(x)=1~\mbox{or}~2$, $\mu(y)=1~\mbox{or}~2$, and $\mu(z)=1~\mbox{or}~2$.
We consider one case, all other cases are similar.

Assume $\mu(x)=\mu(z)=2$ and $\mu(y)=1$. Using the kernel $K$ obtained in Example \ref{6.3.1}, we have
$$
\begin{array}{l}
\Xi_1=I_2\otimes (1,0)\otimes (0,1)=\begin{bmatrix}
0&1&0&0&0&0&0&0\\
0&0&0&1&0&0&0&0\\
\end{bmatrix},\\
\Xi_2=(0,1)\otimes I_2\otimes (0,1)=\begin{bmatrix}
0&0&0&0&0&1&0&0\\
0&0&0&0&0&0&0&1\\
\end{bmatrix},\\
\Xi_3=(0,1)\otimes (1,0) \otimes I_2=\begin{bmatrix}
0&0&0&0&1&0&0&0\\
0&0&0&0&0&1&0&0\\
\end{bmatrix}\\
\end{array}
$$
Then
$$
x=\Xi_1\xi=\begin{bmatrix}
2c_4+2c_5(1-\lambda) +c_6(1-\lambda)\\
-c_4+c_5(\lambda-1)\\
\end{bmatrix}:=\begin{bmatrix}
0\\
1\\
\end{bmatrix}.
$$
$$
y=\Xi_2\xi=\begin{bmatrix}
c_3\\
c_6\lambda
\end{bmatrix}:=\begin{bmatrix}
1\\
\theta\\
\end{bmatrix}.
$$
$$
z=\Xi_3\xi=\begin{bmatrix}
c_2\\
c_3
\end{bmatrix}:=\begin{bmatrix}
0\\
1\\
\end{bmatrix}.
$$
Choosing $\theta$ arbitrary and setting
$$
c_1=c_2=c_4=0;\quad c_3=1; c_6=\theta; c_5=-\frac{1}{2}c_6; \lambda=1,
$$
a direct computation shows that $w=xyz$ is a U-eigenvector of ${\cal A}$ w.r.t. $\lambda=1$ and type ${\cal B}$, where
$$
x=(0,1)^{\mathrm{T}},\quad y=(1,\theta)^{\mathrm{T}},\quad z=(0,1)^{\mathrm{T}}.
$$
\end{exa}

\subsection{Iteration Algorithm}

Only in some simple examples the solution of (\ref{0.0.1}) can easily been obtained from the  solution of (\ref{0.0.2}). We propose the following Iteration Algorithm to solve this problem, which is depicted in Fig. \ref{Fig.6.2}.

\begin{figure}
\centering
\setlength{\unitlength}{0.8 cm}
\begin{picture}(13,19)
\thicklines
\put(0,17){\framebox(6,2){$K(\lambda)=(A-\lambda B)^{\perp}$}}
\put(3,17){\vector(0,-1){2}}
\put(0,13){\framebox(6,2){D-Kernel:$D(\lambda)\subset K(\lambda)$}}
\put(0,10.5){\framebox(2,2){$x(0),\lambda(0)$}}
\put(2.8,11.4) {$\oplus$}
\put(3,13){\vector(0,-1){1.3}}
\put(3,11.3){\vector(0,-1){1.3}}
\put(2,11.5){\vector(1,0){0.8}}
\put(0,8){\framebox(6,2){$\Xi^e_i (x(k))=x_i(k) $}}
\put(3,8){\vector(0,-1){2}}
\put(0,4){\framebox(6,2){$D(\lambda_k)\xi(k)=\prod_{i=1}^rx_i(k)$}}
\put(3,4){\vector(0,-1){2}}
\put(0,0){\framebox(6,2){Lease Square Solution: $\xi(k)$}}
\put(6,1){\vector(1,0){1}}

\put(7,0){\framebox(6,2){$\xi_i(k+1):=\Xi^e_i(D(k)\xi(k))$}}
\put(10,2){\vector(0,1){2}}
\put(7,4){\framebox(6,2){$x(k+1):=\frac{\dsum_{i=1}^r\xi_i(k)}{\|\dsum_{i=1}^r\xi_i(k)\|}$}}
\put(10,6){\vector(0,1){2}}
\put(7,8){\framebox(6,2){$\|x(k+1)-x(k)\|<\epsilon~?$}}
\put(10,10){\vector(0,1){3}}
\put(7,13){\framebox(6,2){Solution: $(\lambda, x)$}}
\put(10,11.5){\line(-1,0){3.5}}
\put(6.5,11.5){\line(0,-1){2.5}}
\put(6.5,9){\vector(-1,0){0.5}}
\put(10.4,11) {Yes}
\put(9.2,11) {No}
\end{picture}
\caption{Iteration Algorithm\label{Fig.6.2}}
\end{figure}

\section{Illustrative Examples}

One of the motivation for proposing the U-eigenproblem is based on the following observation: A hypervector represents a subspace. Hence a U-eigenproblem means for a geven hypermatrix ${\cal A}$, searching a subspace, such that ${\cal A}$ can map it onto a required image subspace. The properties of the image subspace  is determined by the type ${\cal B}$.
We give some examples for describing this.

\begin{exa}\label{e7.1}
\begin{itemize}
\item[(i)] Consider a square matrix $A$. The classical eigenproblem means to find an invariat subspace, which is a invariant. In our U-eigenproblem, to get this purpose, we can simply choose the type $B=I_n$. It is obvious that when our  definition is used to matrix case, it leads to the definition of general eigenvalue-eigenvector.

In this case, one sees that the general eigenvector does not implies an invariant subspace. But depending on the type, it may assign a certain meaning for the image. Say, we may require the image subspace is orthogonal to the original eigenvector.  Then we can choose the type
    \begin{align}\label{7.1}
    B=\begin{bmatrix}0&1\\
    -1&0\end{bmatrix}.
    \end{align}
    Then to find orthogonal eigenvector of $A$, we have the eigenequation as
    \begin{align}\label{7.2}
    Ax=\lambda Bx,
    \end{align}
where $B$ is the one in (\ref{7.1}).

Say,
$$A=\begin{bmatrix}
1&2\\0&1
\end{bmatrix}.
$$
Then we have
$$
(A-\lambda B)x=0.
$$
Equivalently,
$$
(AB^{-1}-\lambda I)B^{-1}x=0.
$$
It is easy to verify that $1$ is a (classical) eigenvalue of $AB^{-1}$ with corresponding eigenvector $z=(1,-1)^{\mathrm{T}}$. Then for original orthogonal eigenproblem, we have $1$ as its eigenvalue with its eigenvector $x=Bz=(-1,1)^{\mathrm{T}}$.

\item[(ii)] Assume ${\cal A}\in \R^{3\times 3\times 3}$. We want to find a vector $x\in \R^3$ such that ${\cal A}$ maps $x$ to  a subspace $S=S(y)$ where $y=y_1y_2$, such that $x\bot S(y)$. We give a way to design this.

Define
$$
B_1=\begin{bmatrix}
0&1&0\\
-1&0&0\\
0&0&1
\end{bmatrix};\quad
B_2=\begin{bmatrix}
0&0&0\\
0&0&1\\
0&-1&0
\end{bmatrix}.
$$
Let
$$
y_i=B_ix,\quad i=1,2,
$$
and define
$$
y=y_1y_2.
$$
Given ${\cal A}\in \R^{3\times 3\times 3}$ with
$a_{1,1,1}=a_{2.2.2}=1$, and other entries being zero.
Its matrix expression $A\in {\cal M}_{9\times 3}$. We consider the eigenproblem as finding an eigenvector $x$ such that its image plan $S(y)$, which is orthogonal to $x$? The problem is formulated as
\begin{align}\label{7.3}
Ax=\lambda y,
\end{align}
where
$$
y=(B_1x)(B_2x)=B_1(I_3\otimes B_2)x^2:=Bx^2.
$$
Then (\ref{7.3}) can be expressed into a standard eigenequation as
\begin{align}\label{7.4}
Ax=\lambda Bx^2.
\end{align}
Furthermore, $x=\Xi^e_{1;3}x^2$, then we have
\begin{align}\label{7.5}
(\tilde{A}-\lambda B)x^2.
\end{align}

where $\tilde{A}=A\Xi^e_{1;3}$ is type depending.

We consider all possible cases.

\begin{itemize}
\item[] Case 1:  Assume the solution type $\mu(x)=1$:
Then $\Xi^1_{(1;3)}=I_3\otimes (1,0,0)$. We have
$\tilde{A}_1=(a_{i,j})$, with $ a_{1,1}=a_{9,7}=1$, and others are zero.
\item[] Case 2:  Assume the solution type $\mu(x)=2$:
Then $\Xi^2_{(1;3)}=I_3\otimes (0,1,0)$. We have
$\tilde{A}_1=(a_{i,j})$, with $a_{1,2}=a_{9,8}=1$, and others are zero.
\item[] Case 3:  Assume the solution type $\mu(x)=3$:
Then $\Xi^3_{(1;3)}=I_3\otimes (0,0,1)$. We have
$\tilde{A}_1=(a_{i,j})$, with $a_{1,3}=a_{9,9}=1$, and others are zero.
\end{itemize}

We consider Case 2, which provides a required solution and the other two Cases do not. (This is because of our choice of $B_1$ and $B_2$.) Using MatLab, it is easy to veriy that (i) $\tilde{A}+\lambda B$ has genetic range $r_g=5$; (ii) the only real essential eigenvalue is $0$; (iii) Corresponding to $0$ the kernel is
$$
K=\Span\{\d_9^1,\d_9^3,\d_9^4,\d_9^5,\d_9^6,\d_9^7,\d_9^9\}.
$$
Finally the only monic D-eigenvector is obtained as
$$
x=\d_3^2.
$$
Then it is easy to calculate that
$$
\begin{array}{l}
y_1=B_1x=\d_3^1\\
y_2=B_2x=-\d_3^3\\
\end{array}
$$
We conclude that $S(y)\perp x$.
\end{itemize}
\end{exa}

Next example comes from \cite{kol14}, and has been solved in \cite{din15}. We use out algorithm to provide a complete solution.

\begin{exa}\label{e7.101} Given
${\cal A}\in \R^{3\times 3\times 3\times 3}$,  and let ${\bf i}=\{1\}$,  ${\bf j}=\{2,3,4\}$. Then
$A=(a_{i,j})=M^{{\bf i}\times {\bf j}}({\cal A})\in {\cal M}_{3\times 27}$ with all nonzero entries as
$$
\begin{array}{llll}
a_{1,1}=0.2883,&a_{1,2}=-0.0031&a_{1,3}=0.1973&a_{1,5}=-0.2485\\
a_{1,6}=-0.2939,&a_{1,9}=0.3847&a_{1,14}=0.2972&a_{1,15}=0.1862\\
a_{1,8}=0.0919,&a_{1,27}=-0.3619&a_{2,14}=0.1241&a_{2,15}=-0.3420\\
a_{2,18}=0.2127,&a_{2,27}=0.2727&a_{3,27}=-0.3054&~\\
\end{array}
$$
$$
Bx^3= x(xTx).
$$

\begin{itemize}
\item[(1)] Find D-eigenvector:
\end{itemize}
The eigenequation is
\begin{align}\label{7.101}
\begin{cases}
{\cal A}x^3=\lambda {\cal B}x^3,\\
 {\cal B}x^3= x(x^{\mathrm{T}}x).
\end{cases}
\end{align}

\begin{rem}\label{r7.101} Note that it is obvious that ${\bf B}x^3=x(x^{\mathrm{T}}x)$ implies ${\bf B}x^4= (x^{\mathrm{T}}x)$.
In addition the norm condition $x^Tx=1$ can be satisfied by finally normalizing the obtained solution $x$.
\end{rem}

First, (\ref{7.101}) can be expressed  into matrix form as
 \begin{align}\label{7.101}
Ax^3=\lambda Bx^3,
\end{align}
where $B\in {\cal M}_{3\times 27}$ with nonzero entries as
$$
b_{1,1}=b_{1,5}=b_{1,9}=b_{2,2}=b_{2,14}=b_{2,18}=b_{3,3}=b_{3,15}=b_{3,27}=1.
$$

To simplify the computation, an observation is: Assume the solution is: $(x_1,x_2,x_3)$. A straightforward computation shows that
$$
Ax^3=(\times,\times,-0.3054x_3^3)^{\mathrm{T}};\quad Bx^3=(\times,\times, x_3^3+x_1^2+x_2^2)^{\mathrm{T}},
$$
where $\times$ stands for uncertain numbers.

w.l.g., we can assume $x_3^3+(x_1^2+x_2^2)\neq 0$.  (Otherwise, we replace $x$ by $-x$.) Then we have
 \begin{align}\label{7.102}
\lambda=\frac{-0.3054x_3^3}{x_3^3+(x_1^2+x_2^2)}.
\end{align}

\begin{itemize}
\item[] Case 1.1 $x_3=0$, $x_1=0$
\end{itemize}

There is no solution.

It is worthy pointing that this eigenequation is homogeneous, so if $x$ is an eigenvector w.r.t. $\lambda$, than $kx$ is also an eigenvector w.r.t. the same $\lambda$. In non-homogeneous case $kx$ remains an eignevector but its corresponding eigenvalue will be $\lambda^s$, where $s$ (might not integer) can also be easily calculated.

\begin{itemize}
\item[] Case 1.2 $x_3=0$, $x_1\neq 0$.
\end{itemize}

Then $\lambda=0$. w.l.g., we assume $x_1=1$. Then we have
$$
A\begin{bmatrix}1\\x_2\\0\end{bmatrix}=0.
$$
There is no solution.

\begin{itemize}
\item[] Case 2.1 $x_3\neq 0$, $x_1=0$.
\end{itemize}

Now $\lambda=\frac{-0.3054}{1+x_2^2}$. Pluging it into (\ref{7.101}) yields two polynomial equations about $x_2$. Numerical computation shows there is no common solution for these to degree 5 polynomial equations.

\begin{itemize}
\item[] Case 2.2 $x_3\neq 0$, $x_1\neq 0$.
\end{itemize}

Assume the eigenvalue is monic, that is, $x_1=0$.

We are looking for only real eigenvectors. Then using (\ref{7.102}) one sees that the corresponding eigenvalue is real.
Note that $A-\lambda B$ has generic rank $3$, we solve the essential eigenvalue by solving
$$
\det[(A-\lambda B)(A-\lambda B)^{\mathrm{T}}]=0.
$$
Calculate it out we have
$$
0.0091\lambda^6-0.0558\lambda^5+0.2397\lambda^4+0.6669\lambda^3+7.6954\lambda^2+8.2062\lambda+27>0.
$$
So there is no essential eigenvalue corresponding to real eigenvector.

Hence we have only to consider quasi-eigenvalues, which corresponding to $K_q$.

It is easy to calculate that the genetic rank of $A-\lambda B$ is $3$, then $\dim(K_q)=24$.

We look for the basis of $K_q$. It is easy to see that
$$
k_i=\d_{27}^i,\quad i=4,7,8,10,11,12,13,16,17,19,20,21,22,23,24,25,26,
$$
are $17$ $\lambda$-independent base element.

The other base elements are $\lambda$-depending. Set
\begin{itemize}
\item[(1)]
$$
Z(i)=
\begin{cases}
A_{1,6},\quad i=1,\\
\lambda*B_{1,1}-A_{1,1}, \quad i=6,\\
0,\quad \mbox{Otherwise}.
\end{cases}
$$
Then set $k_{18}=Z$.
\item[(ii)]
$$
Z(i)=
\begin{cases}
A_{1,6},\quad i=5,\\
\lambda*B_{1,5}-A_{1,5}, \quad i=6,\\
0,\quad \mbox{Otherwise}.
\end{cases}
$$
Then set $k_{19}=Z$.
\item[(iii)]
$$
Z(i)=
\begin{cases}
A_{1,6},\quad i=9,\\
\lambda*B_{1,9}-A_{1,9}, \quad i=6,\\
0,\quad \mbox{Otherwise}.
\end{cases}
$$
Then set $k_{20}=Z$.
\item[(iv)]
$$
C=\begin{bmatrix}
A_{1,2}&A_{1,4}\\
-\lambda B_{2,2}&A_{2.18}-\lambda B_{2,18}\\
\end{bmatrix},
$$
$$
E=[A_{1,18},A_{2,18}-\lambda B_{2,18}],
$$
$$
F=-C^{-1}E,
$$
$$
Z(i)=
\begin{cases}
F(1),\quad i=2,\\
F(2), \quad i=14,\\
0,\quad \mbox{Otherwise}.
\end{cases}
$$
Then set $k_{21}=Z$.
\item[(iv)]
$$
C=\begin{bmatrix}
A_{1,2}&A_{1,6}\\
-\lambda B_{2,2}&0\\
\end{bmatrix},
$$
$$
E=[A_{1,14}A_{2,14},-\lambda B_{2,14}],
$$
$$
F=-C^{-1}E,
$$
$$
Z(i)=
\begin{cases}
F(1),\quad i=1,\\
F(2), \quad i=6,\\
1,\quad i=14,\\
0,\quad \mbox{Otherwise}.
\end{cases}
$$
Then set $k_{22}=Z$.
\item[(v)]
$$
C=\begin{bmatrix}
A_{1,3}&A_{1,6}&A_{1,14}\\
0&0&A_{2,14}-\lambda B_{2,14}\\
-\lambda B_{3,3}&0&0\\
\end{bmatrix},
$$
$$
E=[A_{1,15}, A_{2,15},-\lambda B_{3,15}],
$$
$$
F=-C^{-1}E,
$$
$$
Z(i)=
\begin{cases}
F(1),\quad i=3,\\
F(2), \quad i=6,\\
F(3),\quad i=14,\\
1,\quad i=27,\\
0,\quad \mbox{Otherwise}.
\end{cases}
$$
Then set $k_{23}=Z$.
\item[(vi)]
$$
C=\begin{bmatrix}
A_{1,3}&A_{1,6}&A_{1,14}\\
0&0&A_{2,14}-\lambda B_{2,14}\\
-\lambda B_{3,3}&0&0\\
\end{bmatrix},
$$
$$
E=[A_{1,27}, A_{2,27},A_{3,27}-\lambda B_{3,27}],
$$
$$
F=-C^{-1}E,
$$
$$
Z(i)=
\begin{cases}
F(1),\quad i=3,\\
F(2), \quad i=6,\\
F(3),\quad i=14,\\
1,\quad i=27,\\
0,\quad \mbox{Otherwise}.
\end{cases}
$$
Then set $k_{24}=Z$.
\end{itemize}

Finally, we can find all 24 linearly independent vectors $K_i$ as the basis of $K_q$.
Denote
$$
K=[k_1,\cdots,k_{24}].
$$
They form the basis of candidates of eigenvectors, because $(A-\lambda B)k_i=0$.

Define
$$
W=\begin{bmatrix}
\Xi^e_{(1,3)}-\Xi^e_{(2,3)}\\
\Xi^e_{(1,3)}-\Xi^e_{(3,3)}\\
\end{bmatrix}.
$$
Then $x\in \Col(A-\lambda B)$ has equal projection property. Than is, if $x=x_1x_2x_3$, then$x_1=x_2=x_3$. Since $\rank(W)=4$,
We can choose $20$ elements from $K$, denoted by $N$, such that $(A-\lambda B)N\in W^{\perp}$. (It is an easy linear algebraic computation for finding $N$, we skip the detail here.)

Now the problem becomes finding $\xi\in \R^20$ such that
 \begin{align}\label{7.105}
(A-\lambda B)N\xi=(A-\lambda B)\dsum_{i=1}^{20}\xi_i\Col_i(N)=0.
\end{align}

We already known that there is no exact solution. So we use Iteration Algorithm (Please refer to Figure \ref{Fig.6.2}) to solve this.

A solution is given in \cite{kol14} is $x=(0.5915,-0.7467,-0.3043)^{\mathrm{T}}$. The best corresponding $\lambda$, which makes the error minimum is the least square solution for $A-\lambda B=0$, which is $\lambda=0.1163$.
We use them as the initial eigenvalue and eigenvector, using Iteration Algorithm to find the least square solution.

\begin{table}[!htb]
\centering
\caption{Iteration Result\label{tbhga.2.3a}}
\vskip 2mm
\doublerulesep 0.5pt
\begin{tabular}{|c||c|c|c|}
\hline
$N_0$&Eigenvalue&Eigenvector&Error\\
\hline
\hline
0&-0.1163&(0.5915,-0.7467,-0.3043)&0.1787\\
\hline
1&0.0009&(0.6808,-0.7214,-0.1271)&0.1490\\
\hline
2&0.0002&(0.7625,-0.6424,-0.0770)&0.0572\\
\hline
3&0.0001&(0.7830,-0.6188,-0.0630)&0.0334\\
\hline
4&0.0001&(0.7922,-0.6077,-0.0566)&0.0248\\
\hline
99&4.627e-5&(0.8021,-0.5951,-0.0495)&0.0205\\
\hline
100&4.627e-5&(0.8021,-0.5951,-0.0495)&0.0205\\
\hline
\end{tabular}
\end{table}

It is clear that when the Iteration Algorithm converges, the result is the least square solution of eigenvalue-eigenvector.
We conclude that
$$x=[0.8021,-0.5951,-0.0495)^{\mathrm{T}}]^3$$
is the least square solution of D-eigenvector w.r.t. $\lambda=4.627e-5$.

\begin{itemize}
\item[(2)] Find U-eigenvector:
\end{itemize}

Unlike D-eigenvector case, there are infinite many of U-eigenvector.

Consider the Equation (\ref{7.101}) for non-diagonal case. Then we have
\begin{align}\label{7.111}
\begin{cases}
{\cal A}xyz=\lambda {\cal B}xyz,\\
 {\cal B}xyz= x(y^{\mathrm{T}}z).
\end{cases}
\end{align}
We search the possible solutions. The process is similar to diagonal case. First, we find the kernel of $A-\lambda B$, which is $K$.
Note that
$$
k_j=\d_{27}^{i_j},\quad j\in [1,17],i_j=(4,7,8,10,11,12,13,16,17,19,20,21,22,23,24,25,26).
$$
Then we have
$$
\xi_j=(xyz)_{i_j},\quad j\in [1,17],
$$
that is,
$$
\xi_1=(xyz)_4,\quad \xi_2=(xyz)_7,\quad \cdots,\xi_{17}=(xyz)_{26}.
$$
By deleting first $17$ columns and $i_j$ th rows $i\in [1,17]$, we have remaining equations as
\begin{align}\label{7.112}
K_0\eta=
\begin{bmatrix}
(xyz)_{j_1}\\
(xyz)_{j_2}\\
\vdots\\
(xyz)_{j_{10}}
\end{bmatrix}:=b,
\end{align}
where  $(j_1,j_2,\cdots,j_{10})=(1,2,3,5,6,9,14,15,18,27)$, and
$$
K_0=
\begin{bmatrix}
-0.2939&0&0&0&0&0&0\\
0&0&0&a_{24}&a_{25}&0&0\\
0&0&0&0&0&-1&a_{37}\\
0&-0,2939&0&0&0&0&0\\
a_{51}&a_{52}&a_{53}&0&a_{55}&a_{56}&a_{57}\\
0&0&-0.2939&0&0&0&0\\
0&0&0&a_{64}&1&a_{66}&a_{67}\\
0&0&0&0&0&1&0\\
0&0&0&1&0&0&0\\
0&0&0&0&0&0&1\\
\end{bmatrix},
$$
$a_{ij}$ are $\lambda$-depending entries.
Then we have
\begin{align}\label{7.113}
\begin{cases}
\eta_1=-3.4095 b_1,\\
\eta_2=-3.4095 b_4,\\
\eta_3=-3.4095 b_6,\\
\eta_4=b_9,\\
\eta_5=(1/a_{2,5})b_2-(a_{2,4}/a_{2,5})b_9\\
\eta_6=b_8,\\
\eta_7=b_{10};\\
-b_8+a_{3,7}b_{10}=b_3,\\
-3.4095 a_{5,1}b_1+ a_{5,2}(-3.4095) b_4+a_{5,3}(-3.4095)b_6\\
~~+ a_{5,5}/a_{2,5} b_2-[a_{5,5}a_{2,4}/a_{2,5} +(a_{24}/a_{2,5})]b_9+a_{5,6}b_8+a_{5,7}b_{10}=b_5,\\
[a_{6,4}+(1/a_{2,5})b_2+a_{6,6}b_8+a_{6,7}b_{10}=b_7.
\end{cases}
\end{align}
Next, we need only to choose $x,y,z$ to meet the last three requirements. Note that
$$
\begin{array}{l}
b_1=(xyz)_{j_1}=(xyz)_{1}=x_1y_1z_1,\\
b_2=(xyz)_{2}=x_1y_1z_2,\\
b_3=(xyz)_{3}=x_1y_1z_3,\\
b_4=(xyz)_{5}=x_1y_2z_2,\\
b_5=(xyz)_{6}=x_1y_2z_3,\\
b_6=(xyz)_{9}=x_1y_3z_3,\\
b_7=(xyz)_{14}=x_2y_2z_2,\\
b_8=(xyz)_{15}=x_2y_2z_3,\\
b_9=(xyz)_{18}=x_2y_3z_3,\\
b_{10}=(xyz)_{27}=x_3y_3z_3.\\
\end{array}
$$
\begin{itemize}
\item[] Case 1: Assume $x_1=x_2=0$, $y_3=0$ or $z_3=0$. Then we have $b_i=0$, $i\in[1,10]$, and then $\eta_i=0$, $i\in[1,10]$. It leads to conclusion that
$$
\begin{array}{l}
x=(0,0,1)^{\mathrm{T}};\quad y=(y_1,y_2,0)^{\mathrm{T}};\quad z=(z_1,z_2,z_3)^{\mathrm{T}},\\
x=(0,0,1)^{\mathrm{T}};\quad y=(y_1,y_2,y3)^{\mathrm{T}};\quad z=(z_1,z_2,0)^{\mathrm{T}},
\end{array}
$$
are two sets of (non-diagonal) eigenvectors w.r.t. corresponding $\lambda$.

\item[] Case 2: (general case)
\end{itemize}
Any set of $x,y,z$ satisfying the last three equations of (\ref{7.113}) are eigenvectors.
\end{exa}

Next example is an application to tensors.
\begin{exa}\label{e7.2}

 Let $V$ be a vector space of dimension $n$ with its dual space  $V^*$. A tensor $t$ of covariant order $r$ and contra-variant order $s$ is a multilinear map
	$t:\underbrace{V\times \cdots \times V}_r\times \underbrace{V^*\times \cdots \times V^*}_s\ra \F$. The set of tensors of covariant order $r$ and contra-variant order $s$ is denoted by ${\cal T}^r_s$ \cite{boo86}.

First, we formulate a tensor as a hypermatrix.

Let $\{e_1,\cdots,e_n\}$ be a basis of $V$, and $\{d_1,\cdots,d_n\}$ be its dual basis. Define the structure constants as
$$
t(e_{i_1},\cdots,e_{i_r},d_{j_1},\cdots,d_{j_s})=\mu^{i_1,\cdots,i_r}_{j_1,\cdots,j_s},
		~~i_p,j_q\in [1,n],\quad p\in[1,r],\;q\in [1,s].
$$
Then the set of structure constants form a hypermatrix
\begin{align}\label{7.6}
{\cal A}_t=\left\{\mu^{i_1,\cdots,i_r}_{j_1,\cdots,j_s},
\quad i_p,j_q\in [1,n],\;p\in[1,r],\;q\in [1,s]\right\}
\in \R^{\overbrace{n\times\cdots\times n}^{r+s}}.
\end{align}
Then the structure matrix is
\begin{align}\label{7.7}
\begin{array}{l}
A=M^{{\bf j}\times{\bf i}}({\cal A}_t)\\
=		\begin{bmatrix}
			\mu^{1,\cdots,1}_{1,\cdots,1}&\mu^{1,\cdots,2}_{1,\cdots,1}&\cdots&\mu^{n,\cdots,n}_{1,\cdots,1}\\
			\mu^{1,\cdots,1}_{1,\cdots,2}&\mu^{1,\cdots,2}_{1,\cdots,2}&\cdots&\mu^{n,\cdots,n}_{1,\cdots,2}\\
			\vdots&~&~&~\\
			\mu^{1,\cdots,1}_{n,\cdots,n}&\mu^{1,\cdots,2}_{n,\cdots,n}&\cdots&\mu^{n,\cdots,n}_{n,\cdots,n}\\
		\end{bmatrix}\in {\cal M}_{n^s\times n^r}.
\end{array}
\end{align}
Then it is obvious that
\begin{align}\label{7.8}
t(x_1,\cdots,x_r,\sigma_1,\cdots,\sigma_s)=\ltimes_{j=s}^1 \sigma_j A \ltimes_{i=1}^rx_i.
\end{align}

Next, we consider some geometric structure based on hypermatrix. A square matrix $A$ can be considered as a tensor $t\in {\cal T}^1_1$, with its matrix expression $A$. Then
$$
t(x,\omega)=\omega A x.
$$
$x$ and $\omega$ are  $t$- (or $A$-) orthogonal  if $t(x,\omega)=0$.

This concept can be extended to higher dimension case w.r.t. tensor.
Let $\Omega\subset V^*(\R^n)$ and $D\subset V(\R^n)$ be dual subspace and subspace respectively.
$\Omega$ and $D$ are said to be $t$-orthogonal  if $t(x,\omega)=0$, if for any $\omega_j\in \Omega$, $j\in [1,s]$ and
$x_i\in D$, $i\in [1,r]$, $\omega =\ltimes_{j=s}^i$, $x=\ltimes_{i=1}^rx_i$ satisfies
\begin{align}\label{7.9}
t(x,\omega)=0.
\end{align}
To meet this requirement, we need a type $B_x$ such that
$$
A_{\omega}x=\lambda B x.
$$
And a type $B^*$ such that
$$
\omega A_{x}=\eta \omega B^*.
$$

Then $\Omega$ and $D$ are  $t$-orthogonal, if and only if,
(\ref{7.9}) is true for basis elements, which makes (\ref{7.9}) finitely verifiable.

Particularly, assume $r=s$, then many interesting geometric properties about vectors for matrix case can be extended to subspaces for hypermatrix case.
\end{exa}

Next example shows the potential application in AI.

\begin{exa}\label{e7.3}

Nowadays, transformer obtained wide uses in GAI problems, such as ChatGTP, DeepSeek, etc.  Attention is the engine of a transformer  \cite{vas17}, which is depicted by Fig. \ref{Fig.7.1}.

\begin{figure}
\centering
\setlength{\unitlength}{0.8 cm}
\begin{picture}(8,8)
\thinlines
\put(0,1){\framebox(2,2){}}
\put(3,1){\framebox(2,2){}}
\put(6,1){\framebox(2,2){}}
\put(0,4){\framebox(8,2){}}
\thicklines
\put(0.4,1.2){\vector(0,1){1.6}}
\put(1.6,1.2){\vector(0,1){1.6}}
\put(3.4,1.2){\vector(0,1){1.6}}
\put(4.6,1.2){\vector(0,1){1.6}}
\put(6.4,1.2){\vector(0,1){1.6}}
\put(7.6,1.2){\vector(0,1){1.6}}
\put(0.7,2) {$\cdots$}
\put(3.7,2) {$\cdots$}
\put(6.7,2) {$\cdots$}
\put(3.5,4.8) {${\cal A}$}
\put(7,5.8){\vector(0,1){1.6}}
\put(1,2.8){\vector(0,1){1.6}}
\put(4,2.8){\vector(0,1){1.6}}
\put(7,2.8){\vector(0,1){1.6}}
\put(6,7.4){$V'$}
\put(1,0.2){$Q$}
\put(4,0.2){$K$}
\put(7,0.2){$V$}
\end{picture}
\caption{Attention\label{Fig.7.1}}
\end{figure}

The inputs of an attention are $Q$, $K$, and $V$, which are query, key, and value respectively  \cite{zha23}.

$Q=\{q_1,\cdots,q_r\}$, each $q_i\in \R^n$ is a vector to represent a token, (say, word vector in natural language processing).
Similarly, $K=\{k_1,\cdots,k_r\}$, $V=\{v_1,\cdots,v_r\}$. The attention can be considered as a multi-linear mapping, which is considered as a tensor in AI community. We describe it as a hypermatrix ${\cal A}$.

Then the training process can be described as
\begin{align}\label{7.10}
Aqkv=v',
\end{align}
where $q=\ltimes_{i=1}^rq_i$, $k=\ltimes_{j=1}^rk_j$, and $v=\ltimes_{s=1}^rv_s$, $A$ is the matrix expression of ${\cal A}$.
The optimal training solution is $v'=v$, which means no improvement can be done.

Let $x=qkv$ and ${\cal B}$ be the type, such that $v'={\cal B}(z)$.
And the matrix expression $B$ of ${\cal B}$ can be obtained by using proper $\Xi^e$ (such as we did many times in previous section). Then (\ref{7.10}) can be coverted into the following form.
\begin{align}\label{7.11}
Az=\lambda B z.
\end{align}
So the training becomes to find the eigenvector $z$ for (\ref{7.11}).

In a transformer there are so many hypermatrices, so the technique developed in this paper has potential use there.
\end{exa}

\begin{exa}\label{e7.4} Consider a neural network described in Figure \ref{Fig.7.2}. For simplicity, assume there is only one layer. In fact, multiple layer case is exactly the same.

It is well known that we have
$$
y_j=\dsum_{i=1}^rw_{i,j}x_i,
$$
where $w_{i,j}$ is the weight on the route from $i$ to $j$, and each input $x_i$ and output $y_j$ are vectors.
Assume $w_{i,s,j,t}$ be the coefficient obtained from $x_i^s$ to $y_j^t$, where $x_i^s$ is the $s$-th entry of $x_i$ and $y_j^t$ is the $t$- the entry of $y_j$. Then ${\cal A}=\{w_{i,s,j,t}\}\in \F^{r\times n\times s\times t}$ becomes a hypermatrix .

Then the mapping can be expressed as
\begin{align}\label{7.12}
\begin{cases}
A\ltimes_{i=1}^rx_i=\lambda \ltimes_{j=1}^sy_j,\\
y_j=B_j\ltimes_{i=1}^rx_i,\quad j\in [1,s].
\end{cases}
\end{align}
Where $A$ is the matrix expression of ${\cal A}$, and $B_j$ can be defined to meet our demand. Then to find the required purpose, we seach $x=\ltimes_{i=1}^r x_i$ such that the eigenvector $x$ satisfies
$$
Ax=\lambda Bx.
$$
Comparing  U-Eigenequation (\ref{0.0.1}), one sees that our definition for U-eigenequation is exactly the same as (\ref{7.12}). This fact backups out defenition for U-eigenvalue.

\begin{figure}
\centering
\setlength{\unitlength}{0.8 cm}
\begin{picture}(9,8)
\thinlines
\put(0,2){\framebox(6.5,5.3){}}
\thicklines
\put(1,3){\circle{0.8}}
\put(2.5,3){\circle{0.8}}
\put(3.5,3) {$\cdots$}
\put(5.8,3){\circle{0.8}}
\put(1,6){\circle{0.8}}
\put(2.5,6){\circle{0.8}}
\put(3.5,6) {$\cdots$}
\put(5.8,6){\circle{0.8}}
\put(1,3.4){\vector(0,1){2.25}}
\put(1,3.4){\vector(1,2){1.2}}
\put(1,3.4){\vector(2,1){4.45}}
\put(2.5,3.4){\vector(0,1){2.25}}
\put(2.5,3.4){\vector(-1,2){1.15}}
\put(2.9,3.2){\vector(1,1){2.55}}
\put(5.8,3.45){\vector(0,1){2}}
\put(5.6,3.45){\vector(-2,1){4.2}}
\put(5.6,3.4){\vector(-1,1){2.6}}
\put(1,1){\vector(0,1){1.65}}
\put(2.5,1){\vector(0,1){1.65}}
\put(5.8,1){\vector(0,1){1.65}}
\put(1,6.4){\vector(0,1){1.65}}
\put(2.5,6.4){\vector(0,1){1.65}}
\put(5.8,6.4){\vector(0,1){1.65}}
\put(1,8.2){$y_1$}
\put(1,0.4){$x_1$}
\put(0.1,4.5){$w_{1,1}$}
\put(2.5,8.2){$y_2$}
\put(2.5,0.4){$x_2$}
\put(5.8,8.2){$y_s$}
\put(5.8,0.4){$x_r$}
 \put(5,4.5){$w_{r,s}$}
\put(0.1,6.5){$NN$}
\end{picture}
\caption{A Neural Network\label{Fig.7.2}}
\end{figure}

\end{exa}

\section{Concluding Remarks}

The paper proposed a hypervector form for the eigenvector of hypermatrices. Based on the monic decomposition algorithm, the decomposition of a hypervector becomes finitely computable. Instead of tranditional vector form eigenvector, a hypervector form of eigenvalue represents a subspace. Hence it determines an eigen-subspace. An algorithm has been developed. Comparing with most existing algorithms, which involves the polynomial equations for $x_i^r$, our algorithm  requires to solve polynomial equations for $x_i$ only, that is, multi-linear. Our definition includes the definition for general eigenproblem of matrices  as its special case.

Moreover, the geometric meaning of hypervector eigenvector is obvious, which has been demonstrated in examples.

The followings are some concluding remarks.

\subsection{Hypermatrix vs Tensor}

It is worth to point that hypermatrix and tensor are not the same thing. Only when the basis of a vector space is assigned, the tensor $t\in {\cal T}^r_s$ can be expressed as a hypermatrix. Hence, a hypermatrix can be considered as a tensor only if the basis is pre-assigned. A property of a hypermatrix is called the tensor (-consistent) property, if it is basis-independent.

The eigenvector defined in this paper is a hypervector. Even for the special case, called the D-eigenvector, our definition differ from many existing ones. Say, $x=z^r\in \R^{n^r}$ is an eigenvector, we call $x$ as a (diagonal) eigenvector, while in literature so far, $z\in \R^n$ is called an eigenvector. The difference is, our definition is a tensor property, which is independent on the choice of bases, while the other one is not. So we can use eigenvector to present a subspace, because if $x$ is a eigenvector so is $kx$.

\subsection{Hypermatrix vs Nonsquare Matrix}

Though the eigenproblem considered is about  equilateral hypermatrix, the matrix expression of a hypermatrix is in general non-square. Hence, then eigenproblem considered for hypermatrix is closely related to the eigenproblem for non-square matrix.
For a non-square matrix $A\in {\cal M}_{m\times n}$ a definition for its eigenvalue is \cite{che24}
$$
A\Psi-\lambda I_m=0.
$$
This form is closely related to the eigenequation for hypermatrix. In fact, many results in \cite{che24}, such as generalized Cayley-Hamilton Theorem and invertibility for non-square matrices, etc. can easily be transformed to hypermatrices.

{\bf Acknowledgement}

The authors would like to thank anonymous reviews for their critical comments and helpful suggestions.


\begin{thebibliography}{00}

\bibitem{arz18} D.A. Arzanagh, G. Michailidis, Fast randomized algorithms for t-product based tensor operations and decompositions with applications to imaging data, {\it SIAM J. Imaging Sci.}, Vol. 11, No. 4, 2629-2664, 2018.
%
		\bibitem{bar22}  E. Bar-Shalom, O. Dalin, M. Margaliot, {\it Compound matrices in systems and control theory: a tutorial},  arXiv:2204.00676v1,  2022.
		%
		\bibitem{boo86} W.M. Boothby, {\it Introduction to Differentiable Manifolds and Riemannian Geometry}, 2nd Ed., Elsevier, London, 1986.
%
\bibitem{cche24} C. Chen, {\it Tensor-Based Dynamical Systems, Theory and Applications}, Springer, Switzerland, 2024.
		%
\bibitem{chan22} S. Chang, Y. Wei, T-product tensors, Part II, Tail bounds for sums of random t-product tensors, {\it Comput.   Appl. Math.}, Vol 42, No. 3, 99, 2022.
%
		\bibitem{cha08} K. Chang, K. Oearsib, T. Zhang, {\it Perron-Frobenius theory for nonnegative tensors}, Commun. Math. Sci., 6 (2008), pp. 507-520.
		%
		\bibitem{cha09} K. Chang, K. Oearsib, T. Zhang, {\it On eigenvalue problems of real symmetric tensors}, J. Math. Anal. Appl., 350 (2009), pp. 416-422.
		%
		\bibitem{cha13} K. Chang, L. Qi, T. Zhang, {\it A survey on the spectral theory of nonnegative tensors}, Numer. Linear Algebra Appl., 20 (2013), pp. 891-912.
		%
		\bibitem{che} C. Chen, A.  Surana,  A. M. Bloch, \& I. Rajapakse, {\it Multilinear control systems theory}, SIAM J. Control Optim., 59 (2021), No. 1, pp. 749-776.
		%
		\bibitem{che01} D. Cheng, {\it Semi-tensor product of matrices and its application to Morgen's problem}, Sci. China, Series F, 44 (2001), No. 3, pp. 195-212.
		%
		\bibitem{che07} D. Cheng, H. Qi, {\it Semi-tensor Product of Matrices - Theory and Applications}, Science Press, Beijing, 2007.
		%
		\bibitem{che11} D. Cheng, H. Qi, Z. Li, {\it Analysis and Control of Boolean Networks - A Semi-tensor Product
			Approach}, Springer, London, 2011.
		\bibitem{che12}  D. Cheng, H. Qi, Y. Zhao, {\it An Introduction to Semi-tensor Product of Matrices and Its Applications}, World Scientific, Singapore, 2012.
		%
		\bibitem{che15}  D. Cheng, F. He, H. Qi, T. Xu, {\it Modeling, analysis and control of networked evolutionary games}, IEEE Trans. Autom. Control, 60 (2015), No. 9, pp. 2402-2415.
		%
		\bibitem{che19} D. Cheng, {\it On equivalence of matrices}, Asian J. Math., 23 (2019), No. 2, pp. 257-348.
		%
		\bibitem{che19b} D. Cheng, {\it From Dimension-Free Matrix Theory to Cross-Dimensional Dynamic Systems}, Elsevier, London, 2019.
		%
		\bibitem{che19c} D. Cheng, Z. Liu, {\it A new semi-tensor product of matrices}, Control Theory Technol., 17 (2019), No. 1, pp. 14-22.
%
		\bibitem{che20} D. Cheng, H. Qi, {\it Lectures on Semi-Tensor Product of Matrices, Volume 1: Basic Theory and Multilinear Operation}, Science Press, Beijing, 2020, (in Chinese).
		%
		\bibitem{che21} D. Cheng, Y. Wu, G. Zhao, S. Fu,  {\it A comprehensive survey on STP approach to finite games}, J. Syst. Sci. Complex., 34 (2021), No. 5, pp. 1666-1680.
%
		\bibitem{che21b}  D. Cheng, Z. Ji, J. Feng, S. Fu, J. Zhao, {\it Perfect hypercomplex algebras: Semi-tensor product approach}, Math. Model. Control, 1 (2021), No. 4, pp. 177-187.
		%
		\bibitem{chepr} D. Cheng, Z. Ji, {\it From dimension-free manifolds to dimension-varying control systems}, Comm. Inform. Syst., 23 (2023), No. 1, pp. 85-150.
		%
		\bibitem{che23} D. Cheng,  X. Zhang, Z. Ji, Semi-tensor product of hypermatrices with application to compound hypermatrices,  http:arxiv.org/abs/2303.06295, 2023.
		%
		\bibitem{che23b} D. Cheng, Z. Ji,  {\it Lectures on Semi-Tensor Product of Matrices, Volume 4: Finite and Dimensional-free Dynamic Systems}, Science Press, Beijing, 2023, (in Chinese).
%
		\bibitem{che23c} D. Cheng, M. Meng, X. Zhang, Z. Ji, {\it Contracted product of hypermatrices via STP of matrices},
		Control Theory and Technol., (2023), https://doi.org/10.1007/s11768-023-00155-w.	
		%
\bibitem{che24} D. Cheng, From DK-STP to non-square general Lie algebra and general Lie group, {\it Commun. Inform Sys.}, Vol. 24, No. 1, 1-60, 2024.
%
		\bibitem{chi13} W. Ching, X. Huang, M.K. Ng, t. Siu, {\it Markov Chains: Models, Algorithms and Applications}, 2nd ed., Int. Series in Operations Research and Management Science, 189, Springer, New York, 2013.
		%
		\bibitem{cui14} C. Cui, Y. Dai, J. Nie, {\it All real eigenvalues of symmetric tensors}, SIAM J. Matrix Anal. Appl., 35 (2014), pp. 1582-1601.
		%
\bibitem{del00} L. De Lathauwer, B.D. Moor, J. Vandewalle, A multilinear singular value decomposition, {\it SIAM J. Matrix Anal. Appl.}, Vol. 21, No. 4, 1253-1278, 2000.
%
		\bibitem{din15} W. Ding, Y. Wei, {\it Generalized tensor eigenvalue problems}, SIAM J. Matrix Anal. Appl.,  36 (2015), No. 3, pp. 1073-1099.
		%
		\bibitem{fan20} Z. Fan, C. Deng, H. Li, C. Bu, {\it Multiplications and eignevalues of tensors via linear maps}, Linear Multilinear Algebra, 68 (2020), No. 3, 606-621.
		%
		\bibitem{for16} E. Fornasini, M.E. Valcher, {\it Recent developments in Boolean networks control}, J. Control Decis., 3 (2016), No. 1, pp. 1-18.
		%
\bibitem{goo16} I. Goodfellow, Y. Bengio, A. Courville, {\it Deep Learning}, MIT Press, Massachusetts, 2016.
%
\bibitem{kil13} M.E. Kilmer, K. Braman,  N. Hao, R.C. Hoover, Third-order tensors as operators on matrices: A theoretical and computational framework with applications in imaging, {\it SIAM J. Matrix Anal. Appl.}, Vol. 34, No. 1, 148-172, 2013.
%
\bibitem{kol14} T. G. Kolda, J. R. Mayo, {\it An adaptive shifted power method for computing generalized tensor eigenpairs}, SIAM J. Matrix Anal. Appl., 35 (2014), pp.1563-1581.
	%
		\bibitem{li14} W. Li, M. K. Ng, {\it On the limiting probability distribution of a transition probability tensor}, Linear Multilinear Algebra, 62 (2014), 362-385.
		%
		\bibitem{li18} H. Li, G. Zhao, M. Meng, J. Feng, {\it A survey on applications of semi-tensor product method in engineering}, Sci. China, 61 (2018), 010202:1-010202:17.
		%
		\bibitem{lim05} L. Lim, {\it Singular values and eigenvalues of tensors: A variational approach}, in Proceedings of 1st IEEE International Workshop on Computational Advances in Multi-Sensor Adaptive Processing, 2005, pp. 129-132.
		%
		\bibitem{lim13}
		L. Lim, Tensors and Hypermatrices, in L. Hogben (Ed.) {\it Handbook of Linear Algebra} (2nd ed.), Chapter 15,
		Chapman and Hall/CRC.https://doi.org/10.1201/b16113, 2013.
		%
		\bibitem{lu17}  J. Lu, H. Li, Y. Liu, F. Li, {\it Survey on semi-tensor product method with its applications in logical networks and other finite-valued systems}, IET Control Theory Appl., 11 (2017), No. 13, pp. 2040-2047.
		%
		\bibitem{muh16} A. Muhammad, A. Rushdi, F.A. M. Ghaleb, {\it A tutorial exposition of semi-tensor products of matrices with a stress on their representation of Boolean function}, JKAU Comput. Sci., 5 (2016), pp. 3-30.
		%
		\bibitem{ni14} G. Ni, L. Qi, M. Bai, {\it Geometric measure of entanglement and U-eigenvalues of tensors}, SIAM J. Matrix Anal. Appl., 35 (2014), pp. 73-87.
		%
		\bibitem{pan22} K. Panigrahy, M. Debasisha, Extension of Moore-Penrose inverse of tensor via Einstein product, {\it Linear and Multi-linear Algebra}, Vol. 70, No. 4, 750-773, 2022.
%
		\bibitem{qi05} L. Qi, {\it Eigenvalues of a real supersymmetric tensor},  J. Symbolic Comput., 21 (2005), pp. 1302-1324.
		%
		\bibitem{qi07} L. Qi, {\it Eigenvalues and invariants of tensors}, J. Math. Anal. Appl., 325 (2007), pp. 1363-1377.
		%
		\bibitem{qi08} L. Qi, Y. Wang, E. Wu, {\it D-eigenvalies of diffusion kurtosis tensors}, J. Comput. Appl. Math., 221 (2008), pp. 150-157.
		%
\bibitem{vas17} A. Vaswani, M. Shazeer, N. Parmar, et al., Attention is all you need, {\it Advances in Neural Information Processing Systems},  30, 2017.

%
		\bibitem{wu22} C. Wu, R. Pines, M. Margaliot, J.J. Slotine, {\it Generalization of the multiplicative and additive compounds of square matrices and contraction theory in the hausdorff dimension}, IEEE Trans. Autom. Control, 67 (2022), No. 9, 4629-4644.
%
		\bibitem{yan22} Y. Yan, D. Cheng, J. Feng, H. Li, J. Yue, {\it Survey on applications of algebraic state space theory of logical systems to finite state machines}, Sci. China Inf. Sci., (2022), https://doi.org/10.1007/s11432-22-3538-4.
%
\bibitem{van75} C.F. Van Loan, A general matrix eigenvalue algorithm, {\it SIAM J. Nu,er. Anal.}, Vol. 12, No. 6, 819-834, 1975.

\bibitem{ yan20} H. Yang, Y. Zhao, J. Liu, et la., Sparse regularization tensor robust pca based on t-product and its application in cancer genomic data, {\it Proc. IEEE Int. Conf. Bioinform. Biomedic. (BIBM)}, 2131-2138, 2020.
%
\bibitem{yu22} N. Yu, Z. Liu, R. Gao, Predicting multiple types of microrna-diserse associations based on tensor factorization and label propagation, {\it Computers Biol. Medic.}, Vol. 146, 105558, 2022.
%
\bibitem{zha22} X. Zhang, Z. Ji, D. Cheng, {\it Hidden order of Boolean networks}, {\it IEEE Trans. Neural Networks Learn. Syst.},  https://doi.org/10.1109/TNNLS.2022.3212274, 2022.
%
\bibitem{zj} X. Zhang, M. Meng, and Z. Ji, {\it Analysis of discrete-time switched linear systems under logic dynamic switching},IEEE Transactions on Neural Networks and Learning Systems, 2024. https://ieeexplore.ieee.org/document/10744584
%
\bibitem{zha23} A. Zhang, Z. C. Lipton, A. J. Smola, {\it Dive Into Learning, Pytorch version}, Cambridge University Press, Cambridge, 2023.
%
	\end{thebibliography}
\end{document}